\documentclass[10pt,twocolumn,twoside]{IEEEtran}
\usepackage{amsmath,graphicx}
\usepackage{amssymb}
\usepackage{graphicx}
\usepackage{graphics}
\usepackage{epstopdf}
\usepackage{amsmath}
\usepackage{algorithmicx}
\usepackage{algorithm}
\usepackage{algpseudocode}
\usepackage{multirow}
\usepackage{bm}
\usepackage{enumerate}
\usepackage{paralist}
\usepackage{array}
\usepackage{enumitem}
\usepackage{hhline}

\begin{document}
	
	\title{Approximate Eigenvalue Decompositions of Linear Transformations with a Few Householder Reflectors}
	
	\author{Cristian Rusu\thanks{The author is with the Istituto Italiano di Tecnologia (IIT), Genova, Italy. Contact e-mail address: cristian.rusu@iit.it. Demo source code: https://github.com/cristian-rusu-research/approximate-householder-decomposition}}

	\maketitle
	
	\begin{abstract}
		The ability to decompose a signal in an orthonormal basis (a set of orthogonal components, each normalized to have unit length) using a fast numerical procedure rests at the heart of many signal processing methods and applications. The classic examples are the Fourier and wavelet transforms that enjoy numerically efficient implementations (FFT and FWT, respectively). Unfortunately, orthonormal transformations are in general unstructured, and therefore they do not enjoy low computational complexity properties. In this paper, based on Householder reflectors, we introduce a class of orthonormal matrices that are numerically efficient to manipulate: we control the complexity of matrix-vector multiplications with these matrices using a given parameter. We provide numerical algorithms that approximate any orthonormal or symmetric transform with a new orthonormal or symmetric structure made up of products of a given number of Householder reflectors. We show analyses and numerical evidence to highlight the accuracy of the proposed approximations and provide an application to the case of learning fast Mahalanobis distance metric transformations.
	\end{abstract}
	
	
	\section{Introduction}
	
	The ability to efficiently perform orthonormal linear transformations of data, i.e., complexity $O(n \log n)$ or lower given data of size $n$, is of extreme importance in many practical applications, especially when dealing with high dimensional data or with software running on limited power devices. 
	
	When we discuss computationally efficient linear transformation, the poster algorithm is the fast Fourier transform (FFT) \cite{FFT}. From a computational perspective \cite{FFTComputational}, the FFT is an appropriate way of performing the matrix-vector multiplication between the Fourier matrix $\mathbf{F} \in \mathbb{C}^{n \times n}$ and a given vector. Usually, the matrix-vector multiplication between a general (even orthonormal) matrix and a given vector has quadratic complexity $O(n^2)$ while the FFT uses properties of the highly structured Fourier matrix to reduce the complexity to $O(n \log n)$. The FFT is also related to several other linear transformations that enjoy fast implementations like the discrete cosine transform (DCT) \cite{DCTComputation, DCTComputation2}, the discrete Hartley transform \cite{HartleyComputation}, and the Walsh-Hadamard transform \cite{WalshComputation}.
	
	Another large class of numerically efficient linear transformations is the fast wavelet transforms (FWT) \cite{FWT}. If we consider only orthogonal wavelets \cite{doi:10.1002/cpa.3160410705, 10.1007/978-3-642-97177-8_2}, the first discovered and arguably the most simple and well known is the orthonormal Haar wavelet \cite{HaarCompute}. These algorithms have computational complexity $O(n)$, with an extra speed-up when implemented via a lifting scheme \cite{LiftingScheme}.
	
	Householder reflectors \cite{Householder:1958:UTN:320941.320947} \cite[Chapter~5.1]{Golub1996}\cite{Steinhardt9259} are natural building blocks of orthonormal matrices: an orthonormal transformation of size $n$ has a factorization into $n-1$ reflectors (and a diagonal matrix with entries only in $\{\pm 1 \}$). The idea of building orthonormal transformations that are the product of a given number of Householder reflectors (strictly less than $n-1$) has been studied already in \cite{DictHouseholder} in the context of learning sparsifying transforms. The other basic building blocks of orthonormal matrices have also been considered in similar problems \cite{lemagoarou:hal-01104696, FastSparsifyingTransforms}.

	In this paper, we consider that we are given an orthonormal transformation of size $n$ directly, and our goal is to approximately decompose it into $h$, much smaller than $n$, reflectors. In this case, the proposed solutions are numerically efficient; they do not involve any iterative processes, just an eigenvalue decomposition of the original orthonormal matrix. Next, we use the fact that every symmetric matrix can be diagonalized by orthonormal congruency to propose factorizations of symmetric matrices based on products of a few Householder reflectors, i.e., we factor (approximately) the orthonormal eigenspace of a symmetric matrix by a few Householder reflectors. Unfortunately, the optimization problems that arise in this case (the choice of the $h$ Householder reflectors) are hard and have no closed-form solution. Therefore, we propose iterative optimization algorithms that improve (in the sense that they lower the defined approximation error) the proposed factorization with each step.
	
	It is worth noting that our focus is not on the numerical efficiency of the factorization algorithms for orthonormal or symmetric matrices - current decomposition algorithms exhibit very good computational properties already. Our goal is to construct new factorizations (still, in polynomial time) that approximate well the given matrices and allow for their fast manipulation: for example, after computing the factorization, we want $O(n \log n)$ numerical complexity for matrix-vector multiplication, solving inverse problems, etc.
	
	The paper is organized as follows: Sections II and III deal with the proposed factorizations of orthonormal and symmetric matrices, respectively; Section IV shows synthetic numerical results on random orthonormal and symmetric matrices and finally an application to the construction of fast transforms for distance metric learning.
	
	\section{Approximations of orthonormal matrices}
	
	In this section, we describe two ways of approximating an orthonormal matrix by a product of a fixed, given, number of Householder reflectors together with theoretical insights on the accuracy of such approximations.

	\subsection{The proposed factorization}
	
	In this paper we propose methods to factorize a given orthonormal matrix $\mathbf{U} \in \mathbb{R}^{n \times n}$ into a product of $h$ Householder reflectors such that this factorization is as close as possible to the originally given matrix. The goal is to build an approximate factorization of $\mathbf{U}$ that is highly structured and therefore computationally efficient to use, for example in matrix-vector multiplications. We want to approximate $\mathbf{U}$ by
	\begin{equation}
	\mathbf{\bar{U}} = \mathbf{D} \prod_{k=1}^{h} \mathbf{U}_k = \mathbf{D} \mathbf{U}_h \mathbf{U}_{h-1} \dots \mathbf{U}_1 ,
	\label{eq:thefactorization}
	\end{equation}
	where $\mathbf{D} = \text{diag}(\mathbf{d})$ is a diagonal orthonormal matrix (i.e., $\mathbf{d} \in \{\pm 1 \}^n$) the $h$ factors are Householder reflectors
	\begin{equation}
	\mathbf{U}_k = \mathbf{I} - 2\mathbf{u}_k \mathbf{u}_k^T, \ \| \mathbf{u}_k \|_2 = 1.
	\label{eq:theReflector}
	\end{equation}
	
	Because Householder reflectors are orthonormal, we have that $\mathbf{\bar{U}}$ is also orthonormal. Since every $n \times n$ orthonormal matrix can be written as a product of $n-1$ Householder reflectors (and an orthonormal diagonal), in this paper we are interested in factorizations like \eqref{eq:thefactorization} where $h \ll n$, e.g., $h$ is $O(\log n)$. Because of this restriction on $h$, in general, we are not able to approximate any orthonormal $\mathbf{U}$ (which has $O(n^2)$ degrees of freedom) by the structure $\mathbf{\bar{U}}$ exactly, but a non-zero error will almost always exist. Still, the goal is to reduce this error as much as possible. Because $\det(\mathbf{U}_k) = -1$ then $\det (\prod_{k=1}^h \mathbf{U}_k) = (-1)^h$ and therefore we use the diagonal $\mathbf{D}$ to ensure that the choice of $h$ does not fix the determinant value of $\mathbf{\bar{U}}$.
	
	Matrix-vector multiplication with the matrix $\mathbf{\bar{U}}$ from \eqref{eq:thefactorization} takes $4nh$ operations. Therefore, both the accuracy of approximating $\mathbf{U}$ and the computational complexity of matrix-vector multiplication with $\mathbf{\bar{U}}$ depend on the choice of $h$. We consider the upper bound $h < \frac{n}{2}$ to keep the computational complexity of using $\mathbf{\bar{U}}$ strictly below the $2n^2$ operations needed for the classic, unstructured, matrix-vector multiplication.
	
	\subsection{Constrained Householder reflectors}
	\label{sec:ConstrainedHouseholder}
	
	Given an orthonormal matrix $\mathbf{U}$ and $h$ the number of reflectors, in this section we analyze the following problem:
	\begin{equation}
	\begin{aligned}
	& \underset{\mathbf{u}_k,\ k=1,\dots,h}{\text{minimize}} & & \| \mathbf{U} - \mathbf{\bar{U}}_1 \|_F^2 \\
	& \ \text{subject to} & & \mathbf{\bar{U}}_1 = \prod_{k=1}^{h} \mathbf{U}_k = \mathbf{I} - 2\sum_{k=1}^h \mathbf{u}_k \mathbf{u}_k^T \\
	& & & \| \mathbf{u}_k \|_2 =1,\ \mathbf{u}_k^T \mathbf{u}_j = 0 \text{ for } k \neq j.
	\end{aligned}
	\label{eq:approxWithReflectors}
	\end{equation}
	
	Because the reflector vectors are orthonormal among each other, i.e., $\mathbf{u}_k^T \mathbf{u}_j = 0$, the objective function in \eqref{eq:approxWithReflectors} is
	\begin{equation}
	\| \mathbf{U} - \mathbf{\bar{U}}_1 \|_F^2 
	=  2n - 2\text{tr}(\mathbf{U})+ 2 \sum_{k=1}^h \mathbf{u}_k^T (\mathbf{U} + \mathbf{U}^T) \mathbf{u}_k.
	\label{eq:expandHouseholderConstrained}
	\end{equation}
	
	Therefore, the approximation error depends on the spectral properties of $\mathbf{Z} = \mathbf{U} + \mathbf{U}^T$ which is symmetric and as such has an eigenvalue decomposition $\mathbf{Z} = \mathbf{V}\text{diag}(\mathbf{z})\mathbf{V}^T, \mathbf{z} \in \mathbb{R}^n,$ where $\mathbf{VV}^T = \mathbf{V}^T \mathbf{V} = \mathbf{I}$. Also notice that $\mathbf{\bar{U}}_1$ is always both orthonormal and symmetric. Let us assume without loss of generality that the real-valued eigenvalues $\mathbf{z}$ are sorted in ascending order and that there are $n_-$ negative and $n_+$ positive eigenvalues (and we have that $n_- + n_+ = n$).
	
	\noindent \textbf{Result 1.} Given $\mathbf{U}$, in order to minimize \eqref{eq:approxWithReflectors} with $h \ll n$, the best $\mathbf{\bar{U}}_1$ is composed of $h = n_-$ Householder reflector vectors $\mathbf{u}_k$ which are the eigenvectors corresponding to the negative eigenvalues of $\mathbf{Z}$ and the approximation error in \eqref{eq:expandHouseholderConstrained} is
	\begin{equation}
	\| \mathbf{U} - \mathbf{\bar{U}}_1 \|_F^2 =  2n - 2\text{tr}(\mathbf{U}) + 2\sum_{k=1}^{n_-} z_k
	= 2n - \sum_{k=1}^n |z_k|,
	\label{eq:bestortho}
	\end{equation}
	where we have used that $2\text{tr}(\mathbf{U}) = \text{tr}(\mathbf{U} + \mathbf{U}^T) = \sum_{k=1}^n z_k$. The approximation of $\mathbf{U}$ by $\mathbf{\bar{U}}_1$ is exact when all eigenvalues $z_k \in \{ \pm 2 \}$ and we use $h = n_-$ Householder reflectors, i.e., the number of reflectors is equal to the number of eigenvalues $z_k$ equal to negative two. If $h > n_-$ then we set $\mathbf{u}_k = \mathbf{0}_{n \times 1}$, i.e., $\mathbf{U}_k = \mathbf{I}$, for $k=n_-+1,\dots,h$ because there is no reflector beyond the first $n_-$ that decreases our objective function.
	
	\noindent \textit{Proof.} Let us denote $\mathbf{\tilde{U}} = \begin{bmatrix}
	\mathbf{u}_1 & \dots & \mathbf{u}_h
	\end{bmatrix} \in \mathbb{R}^{n \times h}$ and notice from \eqref{eq:expandHouseholderConstrained} that $\sum_{k=1}^h \mathbf{u}_k^T \mathbf{Z} \mathbf{u}_k = \text{tr}(\mathbf{\tilde{U}}^T \mathbf{Z} \mathbf{\tilde{U}})$ with $\mathbf{\tilde{U}}^T\mathbf{\tilde{U}} = \mathbf{I}$. By the trace minimization Courant-Fischer theorem \cite[Corollary 4.3.39]{MatrixAnalysis} we have that $\underset{\mathbf{\tilde{U}}^T\mathbf{\tilde{U}} = \mathbf{I}}{\min}\ \text{tr}(\mathbf{\tilde{U}}^T \mathbf{Z\tilde{U}} ) = \sum_{k=1}^{n_-} z_k$. $\hfill \blacksquare$
	
	\noindent \textbf{Remark (The 2-norm of large random $\mathbf{Z}$).} Given any orthogonal $\mathbf{U}$ we have that $\| \mathbf{Z} \|_2 \underset{n \rightarrow \infty }{\longrightarrow} 2$ \cite[Section 2.2.2.]{CollinsMale2011}, which means for us that the maximum allowed eigenvalue of $\mathbf{U}$ is achieved for large enough $n$.$\hfill \blacksquare$
	
	\noindent \textbf{Remark (the complex valued case).} Given a unitary matrix $\mathbf{U} \in \mathbb{C}^{n \times n}$ we have the decomposition $\mathbf{U} = \mathbf{T} \text{diag}(\mathbf{\lambda}) \mathbf{T}^H$, $\mathbf{\lambda} \in \mathbb{C}^n$, and its approximation with \eqref{eq:approxWithReflectors} is done with reflectors constructed from the columns of $\mathbf{T}$, denoted $\mathbf{t}_k$. The best performance in this case is achieved as in \eqref{eq:bestortho} and all the reflector vectors are orthonormal to each other, i.e., $\mathbf{t}_k^H \mathbf{t}_j = 0, k \neq j$.$\hfill \blacksquare$

	\subsection{Unconstrained Householder reflectors}
	\label{sec:UnconstrainedHouseholder}
	
	The previously imposed orthogonality condition between the reflector vectors $\mathbf{u}_k$ can be dropped in order to achieve better approximation accuracy. Now we propose to solve the following optimization problem:
	\begin{equation}
	\begin{aligned}
	& \underset{\mathbf{u}_k,\ k=1,\dots,h}{\text{minimize}} & & \| \mathbf{U} - \mathbf{\bar{U}}_2 \|_F^2 \\
	& \ \text{subject to} \ & & \mathbf{\bar{U}}_2 = \prod_{k=1}^{h} \mathbf{U}_k = \prod_{k=1}^h (\mathbf{I} - 2 \mathbf{u}_k \mathbf{u}_k^T) \\
	& & & \| \mathbf{u}_k \|_2 =1.
	\end{aligned}
	\label{eq:approxWithReflectors2}
	\end{equation}
	
	The goal is to update each reflector $\mathbf{U}_k$ separately in order to reduce the value of the objective function. First, some notation. Consider $\mathbf{Z} = \mathbf{U} + \mathbf{U}^T$ and eigenvalue decompositions $\mathbf{Z} = \mathbf{V}\text{diag}(\mathbf{z})\mathbf{V}^T$ and $\mathbf{U} = \mathbf{T} \text{diag}(\mathbf{\lambda}) \mathbf{T}^{H}$ (with $\mathbf{\lambda} \in \mathbb{C}^n$, $\mathbf{T}^H \mathbf{T} = \mathbf{TT}^H = \mathbf{I}$ because $\mathbf{U}$ is a normal matrix)\footnote{We use the \textit{eig} function provided in Matlab$^\text{\textregistered}$ to construct this factorization. This function does not work as described in this paper: it constructs an orthogonal eigenspace only for distinct eigenvalues. If $\mathbf{U}$ has repeated eigenvalues we explicitly orthogonalize their eigenspace using the QR algorithm.} and notice that
	\begin{itemize}
		\item the spectrum of $\mathbf{Z}$ is bounded, i.e., $-2 \leq z_k \leq 2$,
		
		\item the eigenvalues of $\mathbf{U}$, except for $\{\pm 1 \}$, come in complex conjugate pairs $\lambda_{k, k+1} = \alpha_k \pm i\beta_k$ and therefore $\mathbf{Z}$ has two corresponding real eigenvalues $z_{k, k+1} = 2\alpha_k$.
	\end{itemize}
	
	
	To analyze the performance of $\mathbf{\bar{U}}_2$ let us consider a procedure that sequentially initializes $\mathbf{U}_k$ with $k = 1,\dots,h$. To construct an approximation as close as possible to $\mathbf{U}$, for the first reflector $\mathbf{U}_1$ we have to minimize:
	\begin{equation}
	\| \mathbf{U} \! - \mathbf{U}_1 \|_F^2 \! = \! \| \mathbf{U} - \mathbf{I} + \! 2\mathbf{u}_1 \mathbf{u}_1^T \|_F^2 \! = \! 2(n - \text{tr}(\mathbf{U}))\! +\! 2\mathbf{u}_1^T \mathbf{Z} \mathbf{u}_1,
	\end{equation}
	while we add the second reflector $\mathbf{U}_2$ the goal is to minimize:
	\begin{equation}
	\begin{aligned}
	\| \mathbf{U}& -\! \mathbf{U}_2\mathbf{U}_1  \|_F^2 \!=\! \! \| \mathbf{UU}_1\! \! - \! \mathbf{U}_2  \|_F^2 \! = \! \! \| \mathbf{UU}_1\! \!-\! \mathbf{I} + \!2\mathbf{u}_2 \mathbf{u}_2^T  \|_F^2 \\
	& \quad = 2n - 2\text{tr}(\mathbf{UU}_1) + 2\mathbf{u}_2^T (\mathbf{UU}_1 + \mathbf{U}_1\mathbf{U}^T) \mathbf{u}_2.
	\end{aligned}
	\end{equation}
	Assuming $\mathbf{Z}$ has at least one negative eigenvalue, we distinguish now two possibilities:
	\begin{enumerate}
		\item[Case] 1: If the lowest negative eigenvalue of $\mathbf{Z}$ is $-2$, we initialize $\mathbf{u}_1$ to be the eigenvector of this eigenvalue. We have that $\mathbf{Zu}_1 = -2\mathbf{u}_1$ and therefore $\mathbf{Uu}_1 = \mathbf{U}^T\mathbf{u}_1 = -\mathbf{u}_1$. For $\mathbf{u}_2$ we need to check the spectrum of
	\end{enumerate}
	\begin{equation}
	\begin{aligned}
	\mathbf{UU}_1 + \mathbf{U}_1\mathbf{U}^T& =  \mathbf{U}(\mathbf{I}-2\mathbf{u}_1\mathbf{u}_1^T) + (\mathbf{I}-2\mathbf{u}_1\mathbf{u}_1^T)\mathbf{U}^T \\
	= & \mathbf{U} + \mathbf{U}^T - 2\mathbf{Uu}_1 \mathbf{u}_1^T - 2\mathbf{u}_1\mathbf{u}_1^T\mathbf{U}^T \\
	= & \mathbf{Z} + 2 \mathbf{u}_1\mathbf{u}_1^T + 2\mathbf{u}_1\mathbf{u}_1^T = \mathbf{Z} + 4 \mathbf{u}_1 \mathbf{u}_1^T.
	\end{aligned}
	\label{eq:flipeigenvalue}
	\end{equation}
	
	\begin{enumerate}    
		\item[] This is a rank-one update to $\mathbf{Z}$ that flips the $-2$ eigenvalue corresponding to the eigenvector $\mathbf{u}_1$ to $+2$. The other eigenvalues/eigenvectors remain the same.
		
		\item[Case] 2: If the lowest negative eigenvalue of $\mathbf{Z}$ is not $-2$, then we have the duplicate $z_{k, k+1} = 2\alpha_k$ (for $\alpha_k < 0$) where $\lambda_{k, k+1} = \alpha_k \pm i\beta_k$ with $\ |\lambda_{k,k+1}| = 1$ is a complex conjugate pair of eigenvalues of $\mathbf{U}$ with corresponding complex eigenvectors $\mathbf{t}_k$ and $\mathbf{t}_{k+1} = \mathbf{t}_k^*$ that obey $\mathbf{t}_k^H \mathbf{t}_{k+1} = 0$. We choose the reflector vector
		\begin{equation}
		\mathbf{u}_1 = \frac{1}{\sqrt{2}} (\mathbf{t}_k + \mathbf{t}_{k+1}) = \sqrt{2} \Re(\mathbf{t}_k),
		\end{equation}
		for which we achieve the minimum in $\mathbf{u}_1^T \mathbf{Z} \mathbf{u}_1 = 2\alpha_k$ (with $\|\mathbf{u}_1\|_2 = 1$). We now make the argument that we can always construct $\mathbf{u}_2$ such that $\mathbf{u}_2^T(\mathbf{UU}_1 + \mathbf{U}_1\mathbf{U}^T) \mathbf{u}_2 = -2$ which is equivalent to showing that $\mathbf{u}_2^T\mathbf{UU}_1\mathbf{u}_2 = -1$. Consider the reflector vector
	\end{enumerate}
	\begin{equation}
	\mathbf{u}_2 = (\gamma + i\delta)\mathbf{t}_k + (\gamma - i\delta)\mathbf{t}_{k+1}, \gamma^2 + \delta^2 = \frac{1}{2},\ \gamma, \delta \in \mathbb{R}.
	\label{eq:theu2}
	\end{equation}
	\begin{enumerate}
		\item[]    By direct calculation we have that $\mathbf{u}_2^T \mathbf{U} \mathbf{u}_2 = \alpha_k,\ \mathbf{u}_1^T \mathbf{u}_2 = \sqrt{2} \gamma$, $\mathbf{u}_2^T\mathbf{U} \mathbf{u}_1 = \frac{1}{\sqrt{2}}( (\gamma - i\delta)\lambda_k + (\gamma + i\delta)\lambda_{k+1})$. Finally, expanding for \eqref{eq:theu2} we have
	\end{enumerate}
	\begin{equation}
	\begin{aligned}
	\mathbf{u}_2^T\mathbf{UU}_1 \mathbf{u}_2 = & \mathbf{u}_2^2\mathbf{U}\mathbf{u}_2 - 2\mathbf{u}_2\mathbf{U} \mathbf{u}_1 (\mathbf{u}_1^T \mathbf{u}_2) \\
	= & \alpha_k - 2((\gamma - i\delta)\lambda_k + (\gamma + i\delta)\lambda_{k+1}) \gamma.
	\end{aligned}
	\end{equation}
	\begin{enumerate}
		\item [] We set the expression above to $-1$ by choosing $\gamma = -\frac{ \sqrt{1+\alpha_k} }{2},\ \delta = -\frac{\sqrt{1-\alpha_k}}{2}$,
		and we finally have that
		\begin{equation}
		\mathbf{u}_2 = 2(\gamma \Re(\mathbf{t}_k) - \delta\Im(\mathbf{t}_k)).
		\end{equation}
		
		Therefore, the real-valued $\mathbf{u}_1$ does not just minimize the approximation error in the first step but also sets up the problem such that the reduction in the objective function will then be maximal, $-2$, in the second stage.
		
		Because we found a real valued eigenvector $\mathbf{u}_2$ of $\mathbf{U}\mathbf{U}_1$ with eigenvalue $-1$ we now find ourselves in Case 1 and therefore $\mathbf{U}\mathbf{U}_1\mathbf{U}_2 + \mathbf{U}_2\mathbf{U}_1\mathbf{U}^T$ has a new eigenvalue 2, the $-2$ eigenvalue of $\mathbf{UU}_1 + \mathbf{U}_1\mathbf{U}^T$ gets flipped.
		
		For $\zeta = -\frac{\sqrt{1-\alpha_k}}{2}$ and $\eta = \frac{\sqrt{1+\alpha_k}}{2}$ we can also show that there is an eigenvector $\mathbf{v} = 2(\zeta \Re(\mathbf{t}_k) - \eta\Im(\mathbf{t}_k))$ of $\mathbf{UU}_1$ with eigenvalue $1$. Notice, by direct calculation, that $\mathbf{v}^T \mathbf{u}_2 = 0$ and therefore $\mathbf{v}$ is also an eigenvector of $\mathbf{UU}_1\mathbf{U}_2$ with eigenvalue 1. Finally, this means that the spectra of $\mathbf{UU}_1\mathbf{U}_2 + \mathbf{U}_2\mathbf{U}_1 \mathbf{U}^T$ and $\mathbf{Z}$ are identical except that the previous eigenvalues $z_{k, k+1} = 2\alpha_k$ are now both $2$ for some new eigenvectors denoted $\mathbf{v}$ and $\mathbf{u}_2$ that are in the span of $\mathbf{t}_k$ and $\mathbf{t}_{k+1}$.
	\end{enumerate}
	Therefore, given the orthonormal $\mathbf{U}$, to construct the approximation $\mathbf{\bar{U}}_2$, we perform the full eigenvalue decomposition of $\mathbf{U}$, and we follow Case 1 and Case 2 up to the $h$ reflectors or until we exhaust eigenvalues with negative real components.
	
	\noindent \textbf{Result 2.} Given $\mathbf{U}$, in order to minimize \eqref{eq:approxWithReflectors2} with $h \ll n$, the best $\mathbf{\bar{U}}_2$ is composed of $h = n_-$ Householder reflectors and, assuming the eigenvalues $z_k$ of $\mathbf{Z}$ are sorted in ascending order, its approximation error is given by
	\begin{equation}
	\| \mathbf{U} - \mathbf{\bar{U}}_2 \|_F^2  = 2n - 2\text{tr}(\mathbf{\bar{U}}_2^T \mathbf{U}) =  2n_+ - \sum_{k=n_-+1}^n z_k,
	\label{eq:result2}
	\end{equation}
	where we have used $2\text{tr}(\mathbf{\bar{U}}_2^T \mathbf{U}) =  2n_- + \sum_{k=n_-+1}^n z_k$.
	
	The approximation of $\mathbf{U}$ by $\mathbf{\bar{U}}_2$ is exact when all positive eigenvalues $z_k$ are equal to two, and we use $h = n_-$ Householder reflectors, i.e., the number of reflectors is the number of negative eigenvalues $z_k$.
	
	\noindent \textit{Proof.} As previously explained, there is a subspace of dimension $n_-$ where we can find directions such that the objective function \eqref{eq:result2} is reduced. Notice that $\text{tr}(\mathbf{\bar{U}}_2^T \mathbf{U}) =  n_- + \frac{1}{2} \sum_{k=n_-+1}^n z_k = n_- + \sum_{k=n_-+1}^n \alpha_k$, i.e., it is the sum of the spectrum of $\mathbf{U}$ where all negative eigenvalues are replaced by one, the maximum allowed value in the spectrum of $\mathbf{U}$. Therefore the quantity $\text{tr}(\mathbf{\bar{U}}_2^T \mathbf{U})$ is maximized.$\hfill \blacksquare$
	
	For storing the factorizations (either $\mathbf{\bar{U}}_1$ or $\mathbf{\bar{U}}_2$), with floating-point number representations with $Q$ bits, we need $h(n-1)Q + n$ bits, $(n-1)Q$ for each reflector and $n$ for $\mathbf{D}$.
	
	\subsection{Bound on the expected approximation accuracy}
	
	As previously shown, for a particular orthonormal matrix $\mathbf{U}$, the accuracy of the approximation we construct depends on the spectrum of $\mathbf{Z} = \mathbf{U} + \mathbf{U}^T$. For the unconstrained approximation $\mathbf{\bar{U}}_2$ presented in Section \ref{sec:UnconstrainedHouseholder}, notice that the approximation error is different if we consider $\mathbf{U}$ or $-\mathbf{U}$. This difference is caused by the fact that $-\mathbf{U}_k$ is not a Householder reflector if $\mathbf{U}_k$ is. Therefore, depending on the spectrum of $\mathbf{Z}$ we might consider and fix from the beginning $\mathbf{D} = -\mathbf{I}$ in \eqref{eq:thefactorization}, i.e., we consider $-\mathbf{U}$ instead of $\mathbf{U}$ to maximize the summation term in \eqref{eq:result2}. This discussion is mute for the constrained approximation $\mathbf{\bar{U}}_1$ presented in Section \ref{sec:ConstrainedHouseholder} since \eqref{eq:bestortho} depends on the absolute values of the eigenvalues of $\mathbf{Z}$. This is the case because we can write $\mathbf{\bar{U}}_1 = \mathbf{I} - 2 \mathbf{\tilde{U}} \mathbf{\tilde{U}}^T$ where $\mathbf{\tilde{U}} \in \mathbb{R}^{n \times h}$ has orthonormal columns (the $h$ Householder reflector vectors) and we have that $-\mathbf{\bar{U}}_1 = \mathbf{I} - 2 \mathbf{\hat{U}} \mathbf{\hat{U}}^T$ is also a product of Householder reflectors (this time, $n-h$) whose vectors are orthonormal, i.e., $\mathbf{\hat{U}} \in \mathbb{R}^{n \times (n - h)}$ is such that $\begin{bmatrix}
	\mathbf{\tilde{U}} & \mathbf{\hat{U}}
	\end{bmatrix}$ is a full $n \times n$ orthonormal basis.
	
	Furthermore, given that any orthonormal matrix is diagonalized by $n-1$ reflectors, Results 1 and 2 might seem counterintuitive: only $n_{-}$ reflectors are useful (they decrease the objective function value) in the factorization. This is because when we diagonalize with $n-1$ reflectors we do not have an objective function to compute because we know we have enough reflectors to perfectly diagonalize and reach zero approximation error. But if we calculate the objective function value when we diagonalize we notice that this does not monotonically decrease to zero. This explains why our proposed factorization cannot be further improved after constructing $h = n_{-}$ reflectors.
	
	In this section, we present a worse case result on the average performance of the approximation accuracy achieved with $\mathbf{\bar{U}}_2$.
	
	To generate a random orthonormal matrix we build a matrix with i.i.d. entries from the standard Gaussian distribution and then orthogonalize its column by the QR procedure.
	
	\noindent \textbf{Result 3.} Given a random orthonormal $\mathbf{U} \in \mathbb{R}^{n \times n}$ we can always approximate it by $\mathbf{\bar{U}}_2\mathbf{D}$, where $\mathbf{U}_2$ is a product of $h$ Householder reflectors as in \eqref{eq:approxWithReflectors2} and $\mathbf{D}$ is a diagonal matrix with elements $d_{ii} \in \{ \pm 1 \}$ such that
	\begin{equation}
	\mathbb{E}[\| \mathbf{U} - \mathbf{\bar{U}}_2\mathbf{D} \|_F^2] \leq 2(n - h) - \frac{2\sqrt{2}}{\sqrt{\pi}}\sqrt{n-h}.
	\label{eq:result1}
	\end{equation}
	
	\noindent \textit{Proof.} We consider the Householder reflectors $\mathbf{J}_k$ the ones that start the diagonalization process for $\mathbf{U}$, but only for the first $h$ steps of the process \cite[Section~5.2.1]{Golub1996}, i.e., a partial or incomplete QR factorization applied to $\mathbf{U}$, as
	\begin{equation}
	\mathbf{D} \mathbf{J}_h \dots \mathbf{J}_1 \mathbf{U} = \begin{bmatrix}
	\mathbf{I}_{h \times h} & \mathbf{0}_{h \times (n-h)} \\
	\mathbf{0}_{(n-h) \times h} & \mathbf{D}_1\mathbf{\tilde{U}}
	\end{bmatrix},
	\label{eq:ReflectorStructure}
	\end{equation}
	with $\mathbf{D} = \begin{bmatrix}
	\mathbf{I}_{h \times h} & \mathbf{0}_{h \times (n-h)}\\
	\mathbf{0}_{(n-h) \times h} & \mathbf{D}_1
	\end{bmatrix}$ where $\mathbf{\tilde{U}} \in \mathbb{R}^{(n-h) \times (n-h)}$ is orthonormal and $\mathbf{D}_1 \mathbf{\tilde{U}}$ has positive diagonal elements, i.e., $d_{ii}\tilde{u}_{ii} > 0$ for $i=1,\dots,n-h$. Intuitively, each reflector $\mathbf{J}_k$ introduces zeros under the main diagonal on the $k^\text{th}$ column. When the reflectors introduce zeros in an orthonormal matrix (like $\mathbf{U}$, in our case) then we reach the block structure in \eqref{eq:ReflectorStructure}.
	
	The goal is to bring the block structure in \eqref{eq:ReflectorStructure} as close as possible to the identity matrix, i.e., $\mathbf{D} \mathbf{J}_h \dots \mathbf{J}_1 \mathbf{U} \approx \mathbf{I}$. Therefore we want to minimize the quantity
	\begin{equation}
	\| \mathbf{D} \mathbf{J}_h \dots \mathbf{J}_1 \mathbf{U} - \mathbf{I} \|_F^2 = \| \mathbf{U} - \mathbf{J}_1 \dots \mathbf{J}_h \mathbf{D}\|_F^2,
	\end{equation}
	where we used the fact that orthonormal transformations are invariant in the Frobenius norm and that $\mathbf{D}$ and all the reflectors are symmetric, i.e., $\mathbf{J}_k^T = \mathbf{J}_k$. Based on this result, the proposed approximation to $\mathbf{U}$ is $\mathbf{\bar{U}}_2 \mathbf{D}$ where
	\begin{equation}
	\mathbf{\bar{U}}_2 = \mathbf{J}_1 \dots \mathbf{J}_h.
	\label{eq:theU2}
	\end{equation}
	
	Notice that the only non-zero error term comes from the block $\mathbf{D}_1\mathbf{\tilde{U}}$ in \eqref{eq:ReflectorStructure}. To quantify this error, using the expression of $\mathbf{\bar{U}}_2$ from \eqref{eq:theU2}, the Frobenius norm in \eqref{eq:result1} develops to
	\begin{equation*}
	\begin{aligned}
	\| &\mathbf{U} \! - \! \mathbf{\bar{U}}_2\mathbf{D} \|_F^2 \!  = \!  \| \mathbf{U}\! - \! \mathbf{J}_1 \dots \mathbf{J}_h\mathbf{D} \|_F^2 \! =  \! \|  \mathbf{D} \mathbf{J}_h \dots \mathbf{J}_1 \mathbf{U}\! - \! \mathbf{I} \|_F^2  \\
	&=  \|\mathbf{U}\|_F^2 + \| \mathbf{\bar{U}}_2\mathbf{D} \|_F^2 - 2\text{tr}(\mathbf{D}^T \mathbf{\bar{U}}_2^T \mathbf{U} ) \\
	&= 2n - 2\text{tr}(\mathbf{D} \mathbf{J}_h\dots \mathbf{J}_1 \mathbf{U} ) = 2n - 2(h + \text{tr}(\mathbf{D}_1\mathbf{\tilde{U}})) \\
	&= 2(n-h) - 2\text{tr}(\mathbf{D}_1\mathbf{\tilde{U}}),
	\end{aligned}
	\end{equation*}
	where we have used the definition of the Frobenius norm $\|\mathbf{A} \|_F^2 = \text{tr}(\mathbf{A}^T \mathbf{A})$ and the symmetry and orthogonality properties of the reflectors and of the matrix $\mathbf{D}$.
	
	If we consider random orthonormal matrices $\mathbf{U}$ and assume large $n$, then we have that
	\begin{equation}
	\mathbb{E} \left[\text{tr}(\mathbf{D}_1\mathbf{\tilde{U}}) \right] = \mathbb{E} \left[\sum_{i=1}^{n-h} | \tilde{u}_{ii} | \right] = \sqrt{2\pi^{-1} (n - h )}.
	\end{equation}
	The result follows from the fact that the elements of $\mathbf{\tilde{U}}$ can be viewed as Gaussian random variables with zero mean and standard deviation $n^{-\frac{1}{2}}$ (as the columns of $\mathbf{\tilde{U}}$ are normalized in the $\ell_2$ norm, see \cite{OrthoTheory} for an argument on how entries of a random orthonormal matrix weekly converge to the standard Gaussian distribution as $n \to \infty$) and because the $\ell_1$ norm of a standard Gaussian random vector of size $n$ is $\sqrt{2\pi^{-1}}n$. From this observation, the second term on the right hand side of the result in \eqref{eq:result1} follows immediately.
	
	Therefore, the result in \eqref{eq:result1} is achieved with equality when the reflectors are the ones used in the QR decomposition to introduce zeros in the first $h$ columns of $\mathbf{U}$. Given that we construct $\mathbf{\bar{U}}_2$ such that it explicitly minimizes the Frobenius norm in \eqref{eq:approxWithReflectors2} it follows that the approximation is necessarily better and therefore the upper bound in \eqref{eq:result1} holds.
	
	We would like to note here that the reflectors $\mathbf{J}_k = \mathbf{I} - 2\mathbf{j}_k \mathbf{j}_k^T$ used in the QR decomposition to introduce zeros below the $k^\text{th}$ diagonal element have the structure $\mathbf{j}_k = \begin{bmatrix} \mathbf{0}_{(k-1)\times 1} ;& \mathbf{\tilde{j}}_k \end{bmatrix}$, meaning that at most $n-k+1$ entries of $\mathbf{j}_k$ are non-zero. Matrix-vector multiplications between the reflector $\mathbf{J}_k$ and a vector take $4(n-k+1)$ operations. This observation allows for the additional possibility of balancing the computational cost by considering sparse reflector vectors in \eqref{eq:thefactorization}. Unfortunately, in this case, no closed form solution seems to be possible and an iterative optimization problem based on the sparse-PCA approach \cite{SparsePCA} should be considered for each reflector $\mathbf{U}_k$.$\hfill \blacksquare$

	\section{Approximation of symmetric matrices}
	
	Similarly to the previous section, we now describe an algorithm for the approximation of a symmetric matrix by a product of a fixed, given, number of Householder reflectors and a diagonal matrix.
	
	\subsection{The proposed factorization}
	
	Let us now consider fast approximations of symmetric matrices. Given a symmetric matrix $\mathbf{S} \in \mathbb{R}^{n \times n}$ the main result that we will use is its eigenvalue factorization as
	\begin{equation}
	\mathbf{S} = \mathbf{U}\text{diag}(\mathbf{s})\mathbf{U}^T,\ \mathbf{UU}^T = \mathbf{U}^T \mathbf{U} = \mathbf{I}, \ \mathbf{s} \in \mathbb{R}^n,
	\label{eq:symmetric}
	\end{equation}
	where $\mathbf{U}$ stores the orthonormal eigenvectors of $\mathbf{S}$ and where we assume w.l.o.g. that the real-valued entries of $\mathbf{s}$ (the eigenvalues of $\mathbf{S}$) are stored in descending order of their magnitudes. Using the factorization in \eqref{eq:thefactorization}, we now propose an approximation of $\mathbf{S}$ as
	\begin{equation}
	\mathbf{\bar{S}} \! = \! \mathbf{\bar{U}} \text{diag}( \mathbf{\bar{s}}) \mathbf{\bar{U}}^T
	\! \! = \! \mathbf{D} \left( \prod_{k=1}^h \mathbf{U}_k \! \right) \text{diag}(\mathbf{\bar{s}}) \left( \prod_{k=h}^{1} \mathbf{U}_k \! \right) \mathbf{D},
	\label{eq:thesbar}
	\end{equation}
	where all matrices $\mathbf{U}_k$ are Householder reflectors and $\mathbf{\bar{s}} \in \mathbb{R}^n$ now stored the eigenvalues of $\mathbf{\bar{S}}$.
	
	Matrix-vector multiplication with the matrix $\mathbf{\bar{S}}$ from \eqref{eq:thesbar} takes $(8h + 1)n$ operations. We have to consider the bound $h < \frac{n}{4}$ to keep the computational complexity of using $\mathbf{\bar{S}}$ strictly below $2n^2$, the regular computational complexity.

	\subsection{The proposed factorization algorithm}
	
	Given any symmetric matrix $\mathbf{S}$ we want to construct the factorization $\mathbf{\bar{S}}$ as \eqref{eq:thesbar} such that it closely approximates $\mathbf{S}$. There are three components that we can choose in this factorization: i) the spectrum $\mathbf{\bar{s}}$ of the approximation; ii) the number of Householder reflectors $h$ and the values of the reflector vectors $\mathbf{u}_k,\ k=1,\dots,h$ and iii) the orthonormal diagonal matrix $\mathbf{D}$. In this section, we explain how to iteratively and separately choose these components while continuously improving the approximation accuracy. Based on these findings, we propose an algorithm to construct $\mathbf{\bar{S}}$ such that $\|  \mathbf{S} - \mathbf{\bar{S}} \|_F^2$ is reduced.
	
	Our first goal is to choose each Householder reflector $\mathbf{U}_k$ (with all other $h-1$ reflectors fixed) sequentially to minimize
	\begin{equation}
	\| \mathbf{S} - \mathbf{\bar{S}} \|_F^2 \! =  \! \| \mathbf{S} - \mathbf{\bar{U}} \text{diag}( \mathbf{\bar{s}}) \mathbf{\bar{U}}^T \|_F^2  \! = \! \| \mathbf{A}_k  - \mathbf{U}_k \mathbf{B}_k  \mathbf{U}_k \|_F^2,
	\label{eq:symmetricWithHouseholder}
	\end{equation}
	where we have defined the symmetric matrices
	\begin{equation}
	\mathbf{A}_k = \left( \prod_{j={k-1}}^1 \mathbf{U}_j \right) \mathbf{D} \mathbf{S} \mathbf{D} \left( \prod_{j=1}^{k-1} \mathbf{U}_j \right),
	\label{eq:ak}
	\end{equation}
	\begin{equation}
	\mathbf{B}_k = \left( \prod_{j=k+1}^h \mathbf{U}_j \right) \text{diag}(\mathbf{\bar{s}}) \left( \prod_{j=h}^{k+1} \mathbf{U}_j \right).
	\label{eq:bk}
	\end{equation}
	Replacing and developing \eqref{eq:symmetricWithHouseholder} for $\mathbf{U}_k = \mathbf{I} - 2\mathbf{u}_k \mathbf{u}_k^T$ we have
	\begin{equation}
	\begin{aligned}
	\| \mathbf{S} - \mathbf{\bar{S}} & \|_F^2 = \| \mathbf{A}_k  - (\mathbf{I} - 2\mathbf{u}_k \mathbf{u}_k^T) \mathbf{B}_k  (\mathbf{I} - 2\mathbf{u}_k \mathbf{u}_k^T) \|_F^2 \\
	= & \| \mathbf{A}_k \|_F^2 + \| \mathbf{B}_k \|_F^2 - 2\text{tr}(\mathbf{A}_k\mathbf{B}_k) + 4C(\mathbf{u}_k) \\
	= & \| \mathbf{s} \|_2^2 + \| \mathbf{\bar{s}} \|_2^2 - 2\text{tr}(\mathbf{A}_k\mathbf{B}_k) + 4C(\mathbf{u}_k),
	\end{aligned}
	\label{eq:objectivefunctionS3}
	\end{equation}
	where we have denoted
	\begin{equation}
	\begin{aligned}
	C(\mathbf{u}_k) =  \mathbf{u}_k^T(\mathbf{A}_k\mathbf{B}_k + \mathbf{B}_k & \mathbf{A}_k) \mathbf{u}_k \\
	& - 2\text{tr}(\mathbf{A}_k \mathbf{u}_k \mathbf{u}_k^T \mathbf{B}_k \mathbf{u}_k \mathbf{u}_k^T).
	\end{aligned}
	\label{eq:theRk}
	\end{equation}
	Finding the $\mathbf{u}_k$ with $\| \mathbf{u}_k \|_2 = 1$ that minimizes \eqref{eq:objectivefunctionS3} seems hard in general (making no assumptions on the spectra of $\mathbf{A}_k$ and $\mathbf{B}_k$) as it seems that there is no closed form solution.
	
	But notice that we can separately optimize the two parts of the expression in \eqref{eq:theRk}:
	\begin{itemize}[leftmargin=*]
		\item Notice that the minimizer of $\mathbf{u}_k^T(\mathbf{A}_k \mathbf{B}_k + \mathbf{B}_k \mathbf{A}_k) \mathbf{u}_k$ is the eigenvector of the smallest, negative, eigenvalue of the symmetric matrix $\mathbf{A}_k \mathbf{B}_k + \mathbf{B}_k \mathbf{A}_k$. We denote it $\mathbf{u}_k^\dagger$.
		
		\item To maximize the second trace term in $\mathbf{u}_k$ we develop:
	\end{itemize}
	\begin{equation}
	\begin{aligned}
	\text{tr}(\mathbf{A}_k \mathbf{u}_k \mathbf{u}_k^T & \mathbf{B}_k \mathbf{u}_k \mathbf{u}_k^T) = \text{tr}(\mathbf{u}_k \mathbf{u}_k^T \mathbf{A}_k \mathbf{u}_k \mathbf{u}_k^T \mathbf{B}_k) \\
	= & \text{vec}(\mathbf{u}_k \mathbf{u}_k^T)^T \text{vec}(\mathbf{A}_k \mathbf{u}_k \mathbf{u}_k^T \mathbf{B}_k ) \\
	= & \text{vec}(\mathbf{u}_k \mathbf{u}_k^T)^T (\mathbf{B}_k \otimes \mathbf{A}_k ) \text{vec}(\mathbf{u}_k \mathbf{u}_k^T) \\
	= & (\mathbf{u}_k \otimes \mathbf{u}_k)^T (\mathbf{B}_k \otimes \mathbf{A}_k ) (\mathbf{u}_k \otimes \mathbf{u}_k).
	\end{aligned}
	\end{equation}
	\begin{itemize}[leftmargin=*]
		\item[] We have used the cyclic property $\text{tr}(\mathbf{XYZ}) = \text{tr}(\mathbf{ZXY})$ and vectorize property $\text{tr}(\mathbf{X}^T \mathbf{Y}) = \text{vec}(\mathbf{X})^T \text{vec}(\mathbf{Y})$ of the trace, the fact that $\text{vec}(\mathbf{XYZ}) = (\mathbf{Z}^T \otimes \mathbf{X})\text{vec}(\mathbf{Y})$ and $\text{vec}(\mathbf{xx}^T) = \mathbf{x} \otimes \mathbf{x}$. We denote $\mathbf{v} \in \mathbb{R}^{n^2}$ the eigenvector corresponding to the highest, positive, eigenvalue of the symmetric matrix $\mathbf{B}_k \otimes \mathbf{A}_k$. We note that the Kronecker product is never explicitly calculated in order to obtain $\mathbf{v}$ but we compute: $\lambda_{\mathbf{A}_k}^{\min}, \lambda_{\mathbf{B}_k}^{\min}$ and $\lambda_{\mathbf{A}_k}^{\max}, \lambda_{\mathbf{B}_k}^{\max}$ the pairs of lowest and highest eigenvalues of $\mathbf{A}_k$ and $\mathbf{B}_k$, respectively and then the highest eigenvalue of $\mathbf{B}_k \otimes \mathbf{A}_k$ is $\max\{ \lambda_{\mathbf{A}_k}^{\min}\lambda_{\mathbf{B}_k}^{\min},\ \lambda_{\mathbf{A}_k}^{\min}\lambda_{\mathbf{B}_k}^{\max},\ \lambda_{\mathbf{A}_k}^{\max}\lambda_{\mathbf{B}_k}^{\min},\ \lambda_{\mathbf{A}_k}^{\max}\lambda_{\mathbf{B}_k}^{\max}  \}$ and the corresponding eigenvector is the Kronecker product of the eigenvectors corresponding to the eigenvalues whose product is maximum. For example, for positive semidefinite matrices $\mathbf{A}_k$ and $\mathbf{B}_k$ the maximum eigenvalue of $\mathbf{B}_k \otimes \mathbf{A}_k$ is $\lambda_{\mathbf{A}_k}^{\max}\lambda_{\mathbf{B}_k}^{\max}$ and therefore its corresponding eigenvector is $\mathbf{v} = \mathbf{v}_{\mathbf{B}_k}^{\max} \otimes \mathbf{v}_{\mathbf{A}_k}^{\max}$, i.e., the Kronecker product of the eigenvectors corresponding to the highest eigenvalues \cite[Chapter~2.1]{OnTheKron}. With the $\mathbf{v}$ just computed we now minimize
	\end{itemize}
	\begin{equation}
	\| \mathbf{v} - \mathbf{u}_k \otimes \mathbf{u}_k \|_F^2 = \| \mathbf{v} - \text{vec}(\mathbf{u}_k\mathbf{u}_k^T) \|_F^2 = \| \mathbf{V} - \mathbf{u}_k \mathbf{u}_k^T \|_F^2,
	\label{eq:factor}
	\end{equation}
	\begin{itemize}[leftmargin=*]
		\item[] where $\mathbf{V}$ is an $n \times n$ matrix build from the re-arranged elements of $\mathbf{v}$ (for details on the arrangement procedure please see \cite[Sections 6 and 7]{ubiquitousKronecker}). Notice that because $\mathbf{v}$ is a Kronecker product, we have that $\mathbf{V}$ is a rank one matrix. The minimizer of \eqref{eq:factor} is the best rank one approximation of $\mathbf{V} + \mathbf{V}^T$ \cite[Section 7]{ubiquitousKronecker}. We denote this solution $\mathbf{u}_k^\ddagger$.
		
		\noindent \textbf{Remark (initialization of $\mathbf{u}_k^\ddagger$ based on generalized Rayleigh quotient calculations when $\mathbf{A}_k$ and $\mathbf{B}_k$ are positive definite).} The second term of $C(\mathbf{x})$ can be written
		\begin{equation}
		R(\mathbf{A}, \mathbf{B}, \mathbf{x}) = \frac{\mathbf{x}^T \mathbf{Ax} \mathbf{x}^T \mathbf{Bx}}{\mathbf{x}^T \mathbf{x}}.
		\end{equation}
		Without loss of generality, for convenience, we momentarily drop the index $k$ from the notation. Assuming the matrices are positive semidefinite and $\mathbf{B}$ is invertible, we denote the Cholesky factorization $\mathbf{B} = \mathbf{LL}^T$, we make the change of variable $\mathbf{y} = \mathbf{L}^T\mathbf{x}$ and we have
		\begin{equation}
		\! R(\mathbf{A}, \mathbf{B}, \mathbf{x}) \! \! = \! \!  \frac{\mathbf{x}^T \mathbf{Ax} \mathbf{x}^T \mathbf{LL}^T\mathbf{x}}{\mathbf{x}^T \mathbf{x}} \! \! = \! \! \frac{ \mathbf{y}^T \mathbf{L}^{-1}\mathbf{A}\mathbf{L}^{-T}\mathbf{y} \mathbf{y}^T \mathbf{y} }{\mathbf{y}^T \mathbf{L}^{-1} \mathbf{L}^{-T} \mathbf{y}}.
		\end{equation}
		This expression is almost a generalized Rayleigh quotient \cite[Chapter~8.2.3]{Golub1996} (we have an extra multiplicative term $\mathbf{y}^T \mathbf{y}$). Assuming $\| \mathbf{y}\|_2 = 1$, to maximize the last quantity in $\mathbf{y}$ we use the generalized eigenvalue decomposition which reduces to finding the eigenvector $\mathbf{y}$ corresponding to the highest eigenvalue of 
		$\mathbf{L}^T \mathbf{A L}^{-T}$ \cite[Chapter~7.7]{Golub1996}. We recover $\mathbf{x} = \mathbf{L}^{-T}\mathbf{y}$ and normalize $\mathbf{u}_k^\ddagger = \| \mathbf{x} \|_2^{-1} \mathbf{x}$.$\hfill \blacksquare$
	\end{itemize}
	Finally, given the two vectors $\mathbf{u}_k^\dagger$ and $\mathbf{u}_k^\ddagger$, we initialize the reflector vector $\mathbf{u}_k^{(1)}$ by a two step procedure. First, we update $\mathbf{u}_k^\ddagger$ such that $(\mathbf{u}_k^\dagger)^T \mathbf{u}_k^\ddagger = 0$ is guaranteed
	\begin{equation}
	\mathbf{u}_k^\ddagger \leftarrow  \mathbf{u}_k^\ddagger - (\mathbf{u}_k^\dagger)^T \mathbf{u}_k^\ddagger \mathbf{u}_k^\dagger,\ \mathbf{u}_k^\ddagger \leftarrow \mathbf{u}_k^\ddagger \| \mathbf{u}_k^\ddagger \|_2^{-1},
	\label{eq:processtheinitial}
	\end{equation}
	and then we solve the minimization problem 
	\begin{equation}
	\mathbf{u}_k^{(1)} = \underset{\mathbf{u}_k = (1-\gamma_k^2/2) \mathbf{u}_k^\dagger \pm \sqrt{(\gamma_k^2 - \gamma_k^4/4)} \mathbf{u}_k^\ddagger}{\arg \min} \quad  C(\mathbf{u}_k),
	\label{eq:combine}
	\end{equation}
	where $\gamma_k \in [0,\sqrt{2}]$ is chosen such that $\mathbf{u}_k^{(1)}$ minimizes $C(\mathbf{u}_k^{(1)})$, i.e., we sweep over the unit hypersphere between points $\mathbf{u}_k^\dagger$ and $\mathbf{u}_k^\ddagger$ in order to minimize $C(\mathbf{u}_k)$. Since we optimize separately the two terms in \eqref{eq:theRk} we search over linear combinations of these for a good minimizer of the overall expression. The formulation of $\mathbf{u}_k$ in \eqref{eq:combine} guarantees that $\| \mathbf{u}_k^{(1)} \|_2 = 1$ for any $\gamma_k$ while the one dimensional minimization (in $\gamma_k$) is done efficiently with a numerical procedure\footnote{We use the \textit{fminbnd} function, initialized at $\gamma_k = 0$, that is provided in Matlab$^\text{\textregistered}$ to find a local minimizer of $C(\mathbf{u}_k)$ as a function of $\gamma_k$.}.
	
	With the initialized reflector vector $\mathbf{u}_k^{(1)}$ we start now an iterative gradient descent procedure to further reduce $C(\mathbf{u}_k)$. Notice that, by the cyclic property of the trace we have $\text{tr}(\mathbf{A}_k \mathbf{u}_k \mathbf{u}_k^T \mathbf{B}_k \mathbf{u}_k \mathbf{u}_k^T) = \text{tr}(\mathbf{u}_k^T\mathbf{A}_k \mathbf{u}_k \mathbf{u}_k^T \mathbf{B}_k \mathbf{u}_k) = (\mathbf{u}_k^T\mathbf{A}_k \mathbf{u}_k) (\mathbf{u}_k^T \mathbf{B}_k \mathbf{u}_k)$ is the product of two quadratic forms (the convexity of such products for optimization purposes was studied in \cite{ProdQuadForms}). As such, we have the gradient expression
	\begin{equation}
	\begin{aligned}
	\nabla C&(\mathbf{u}_k^{(i)}) = 2(\mathbf{A}_k\mathbf{B}_k + \mathbf{B}_k \mathbf{A}_k) \mathbf{u}_k^{(i)} \\
	& - 4 ( (\mathbf{u}_k^{(i)})^T \mathbf{A}_k \mathbf{u}_k^{(i)} \mathbf{B}_k + (\mathbf{u}_k^{(i)})^T \mathbf{B}_k \mathbf{u}_k^{(i)} \mathbf{A}_k )\mathbf{u}_k^{(i)}.
	\end{aligned}
	\label{eq:gradient}
	\end{equation}
	
	We update the gradient in the same fashion as \eqref{eq:processtheinitial} to obtain
	\begin{equation*}
	\mathbf{g}_k^{(i)} = \nabla C(\mathbf{u}_k^{(i)}) -(\mathbf{u}_k^{(i)})^T \nabla C(\mathbf{u}_k^{(i)}) \mathbf{u}_k^{(i)},\ \mathbf{g}_k^{(i)} \leftarrow \mathbf{g}_k^{(i)} \| \mathbf{g}_k^{(i)} \|_2^{-1},
	\end{equation*}
	and finally we have the update equation:
	\begin{equation}
	\mathbf{u}_k^{(i+1)} = (1-\gamma_{\star,k}^2/2)\mathbf{u}_k^{(i)} - \sqrt{(\gamma_{\star,k}^2 - \gamma_{\star,k}^4/4)} \mathbf{g}_k^{(i)},
	\label{eq:iterate}
	\end{equation}
	where $\gamma_{\star,k}$ is found by a one dimensional search such that $C(\mathbf{u}_k^{(i+1)})$ is minimized. Since the search for $\mathbf{u}_k^{(i+1)}$ contains the previous solution $\mathbf{u}_k^{(i)}$, i.e., the step cannot increase the objective function, this iterative procedure has a strictly monotonically descent to a stationary point.
	
	Finally, to update the diagonal of $\mathbf{D}$ in \eqref{eq:thesbar} denoted by $\mathbf{d}$ we minimize $\|  \mathbf{S} - \mathbf{D} \mathbf{B}_0 \mathbf{D} \|_F^2$. If we denote by $\mathbf{\tilde{s}}_i,\ \mathbf{\tilde{b}}_i \in \mathbb{R}^{n-1}$ the $i^\text{th}$ rows of $\mathbf{S}$ and $\mathbf{B}_0$, respectively, both with the diagonal element removed then
	\begin{equation}
	d_{ii} = 1 \text{ if } \| \mathbf{\tilde{s}}_i - \mathbf{\tilde{b}}_i  \|_2 \geq \| \mathbf{\tilde{s}}_i + \mathbf{\tilde{b}}_i \|_2 \text{ else } d_{ii} = -\!1.
	\label{eq:optimaldiagonal}
	\end{equation}
	
	For completeness, the full proposed learning procedure, which we call Symmetric Householder Factorization (SHF), is presented in Algorithm 1. The algorithm runs for fixed $K$ iterations or until the progress in relative error between two consecutive iterations is below $\epsilon = 10e-8$. In our discussion so far we have assumed that the spectra of the given matrix $\mathbf{S}$ and its approximation $\mathbf{\bar{S}}$ are identical. We can always also optimize over the choice of the spectrum $\mathbf{\bar{s}}$ by minimizing
	\begin{equation}
	\| \mathbf{S} - \mathbf{\bar{S}} \|_F^2 \! = \! \|\mathbf{\bar{U}}^T  \mathbf{S} \mathbf{\bar{U}} -  \text{diag}( \mathbf{\bar{s}}) \|_F^2,
	\label{eq:updateofsbar}
	\end{equation}
	which is given by $\mathbf{\bar{s}} = \text{diag}(\mathbf{\bar{U}}^T  \mathbf{S} \mathbf{\bar{U}})$. We can trivially adapt Algorithm 1 to also perform this update iteratively after the calculation of all the reflectors. We call this approach SHF with Spectrum Update (SHF--SU).

	\begin{algorithm}[t]
		\caption{ \textbf{\!--\! Symmetric Householder Factorization (SHF)} \newline \textbf{Input: } The symmetric matrix $\mathbf{S} \in \mathbb{R}^{n \times n}$, the number of Householder reflectors $h$ in the factorization and the maximum number of iterations $K$. \newline \textbf{Output: }  The symmetric matrix $\mathbf{\bar{S}}$ factored as \eqref{eq:thesbar} such that $\| \mathbf{S} - \mathbf{\bar{S}} \|_F^2$ is reduced.}
		\begin{algorithmic}
			\State \textbf{1. } Construct the eigenvalue decomposition of $\mathbf{S}$ as in \eqref{eq:symmetric}.
			
			\State \textbf{2. } Initialize all elements of the approximate factorization:
			
			\begin{itemize}
				\item The spectrum is set to $\mathbf{\bar{s}} = \mathbf{s}$.
				
				\item Set all reflector vectors $\mathbf{u}_k = \mathbf{0},\ k = 1,\dots,h$.
				
				\item Update all $h$ reflector vectors according to \eqref{eq:combine}.
				
				\item The diagonal elements of $\mathbf{D}$ are set according to \eqref{eq:optimaldiagonal}.
			\end{itemize}
			
			\State \textbf{3. } For $1,\dots,K:$
			\begin{itemize}
				\item Iteratively update all reflectors, for $k=1,\dots,h$:
				\begin{itemize}
					\item Construct $\mathbf{A}_k$ and $\mathbf{B}_k$ according to \eqref{eq:ak} and \eqref{eq:bk}.
					
					\item Starting from the current $\mathbf{u}_k$, iteratively update the reflector vector according to \eqref{eq:iterate} until convergence.
				\end{itemize}
				
				\item Update $\mathbf{D}$ according to \eqref{eq:optimaldiagonal}.
			\end{itemize}
		\end{algorithmic}
	\end{algorithm}
	\noindent \textbf{Remark (bounding $C(\mathbf{u}_k)$).} Denoting $\mathbf{y} = \left( \begin{bmatrix} \mathbf{u}_k \\ 1 \end{bmatrix} \otimes \mathbf{u}_k \right) \in \mathbb{R}^{n^2 + n}$, an alternative way of writing the total cost is
	\begin{equation}
	C(\mathbf{u}_k) = \mathbf{y}^T \begin{bmatrix}
	-2(\mathbf{B}_k \otimes \mathbf{A}_k)& \mathbf{0}_{n^2 \times n} \\
	\mathbf{0}_{n \times n^2} & \mathbf{A}_k\mathbf{B}_k + \mathbf{B}_k\mathbf{A}_k 
	\end{bmatrix} \mathbf{y}.
	\end{equation}
	Therefore, if the block diagonal matrix is positive semidefinite (meaning that both $-(\mathbf{B}_k \otimes \mathbf{A}_k)$ and $\mathbf{A}_k\mathbf{B}_k + \mathbf{B}_k\mathbf{A}_k $ are positive semidefinite) then there is no $\mathbf{u}_k$ such that $C(\mathbf{u}_k) < 0$ and therefore the objective function in \eqref{eq:objectivefunctionS3} cannot be reduced (this is a sufficient condition).
	
	Assume that the given $\mathbf{S}$ is positive semidefinite. Notice that in our case the spectra of $\mathbf{A}_k$ and $\mathbf{B}_k$ are identical to the spectrum of $\mathbf{S}$ (since $\mathbf{\bar{s}} = \mathbf{s}$) and then general results from linear algebra \cite[Chapter 4]{MatrixAnalysis} show that the highest eigenvalues of $\mathbf{A}_k \mathbf{B}_k + \mathbf{B}_k \mathbf{A}_k$ and $\mathbf{B}_k \otimes \mathbf{A}_k$ are bounded by $2 (\lambda_{\mathbf{S}}^{\min})^2 \leq \lambda_{\mathbf{A}_k \mathbf{B}_k + \mathbf{B}_k \mathbf{A}_k}^{\max} \leq 2 (\lambda_{\mathbf{S}}^{\max})^2$ and $(\lambda_{\mathbf{S}}^{\min})^2 \leq \lambda_{\mathbf{B}_k \otimes \mathbf{A}_k}^{\max} \leq (\lambda_{\mathbf{S}}^{\max})^2$, respectively, which leads to the bound
	\begin{equation}
	\frac{C(\mathbf{u}_k)}{2( (\lambda_{\mathbf{S}}^{\max})^2 - (\lambda_{\mathbf{S}}^{\min})^2)} \in [-1,1].
	\end{equation}
	When $\mathbf{A}_k$ and $\mathbf{B}_k$ share the same eigenspace we trivially have that $C(\mathbf{u}_k) = 0$ for any $\mathbf{u}_k$ -- this reflects the situation where the $h$ reflectors exactly describe the eigenspace of $\mathbf{S}$.$\hfill \blacksquare$
	
	\noindent\textbf{Remark (analysis of a relaxed problem when $\mathbf{A}_k$ and $\mathbf{B}_k$ are positive definite).} Consider a similarly constrained optimization problem
	\begin{equation}
	\underset{\mathbf{u}_k}{\text{minimize}} \quad C(\mathbf{u}_k) \quad    \text{subject to}  \quad \| \mathbf{u}_k \|_2^2 \leq 1,
	\label{eq:relaxedsolution}
	\end{equation}
	The constraint is important to ensure a bounded solution, otherwise $C(\mathbf{u}_k) = -\infty$. We now defined the Lagrangian
	\begin{equation}
	\mathcal{L}(\mathbf{u}_k, \nu) = C(\mathbf{u}_k) + \nu (\mathbf{u}_k^T \mathbf{u}_k - 1),
	\end{equation}
	with $\nu \in \mathbb{R}_+$, a Lagrange multiplier, and define the derivative
	\begin{equation}
	\nabla \mathcal{L}(\mathbf{u}_k, \nu) = \nabla C(\mathbf{u}_k) + 2\nu \mathbf{u}_k.
	\end{equation}
	If $\mathbf{u}_k^\star$ minimizes $\mathcal{L}(\mathbf{u}_k, \nu^\star)$ then $\nabla \mathcal{L}(\mathbf{u}_k^\star, \nu^\star) = \mathbf{0}$ (the stationarity condition) and therefore $(\mathbf{u}_k^\star)^T \nabla \mathcal{L}(\mathbf{u}_k^\star, \nu^\star) = 0$. Plugin the gradient from \eqref{eq:gradient} leads by direct calculation to
	\begin{equation}
	C(\mathbf{u}_k^\star) = 2(\mathbf{u}_k^\star)^T \mathbf{A}_k \mathbf{u}_k^\star(\mathbf{u}_k^\star)^T \mathbf{B}_k \mathbf{u}_k^\star - \nu^\star \| \mathbf{u}_k^\star \|_2^2.
	\end{equation}
	To complete the KKT conditions for \eqref{eq:relaxedsolution}, we also have that $\nu^\star \geq 0$ (dual feasibility) and $\nu^\star (\|\mathbf{u}_k^\star\|_2^2 - 1) = 0$ (complementary slackness). The optimal solution of \eqref{eq:relaxedsolution} is such that $C(\mathbf{u}_k^\star) \leq 0$, since the feasible $\mathbf{u}_k^\star = \mathbf{0}$ trivially leads to $C(\mathbf{u}_k^\star) = 0$, and therefore $\nu^\star \geq \frac{2(\mathbf{u}_k^\star)^T \mathbf{A}_k \mathbf{u}_k^\star(\mathbf{u}_k^\star)^T \mathbf{B}_k \mathbf{u}_k^\star}{\|\mathbf{u}_k^\star\|_2^2}$. In the special case where $\mathbf{A}_k$ and $\mathbf{B}_k$ are positive definite then we necessarily have that $\nu^\star > 0$ and therefore $\| \mathbf{u}_k^\star \|_2^2 = 1$.
	$\hfill \blacksquare$
	
	\noindent \textbf{Remark (a matrix manifold optimization approach).} The problem at hand can be seen as the minimization of $C(\mathbf{x})$ \eqref{eq:theRk} over the unit sphere manifold $S^{n-1}$ for which iterative optimization procedures are available \cite[Chapters 4.6 and 6.4]{10.2307/j.ctt7smmk} that could be adapted to our case. $\hfill \blacksquare$
	
	\subsection{Bound on the average performance of SHF}
	
	For a particular symmetric matrix $\mathbf{S}$, the accuracy of the approximation we construct depends on its spectrum. In this section, we present a worse case result on the average performance of the approximation accuracy achieved with $\mathbf{\bar{S}}$.
	
	\noindent \textbf{Result 4.} We generate a random symmetric matrix we build a matrix $\mathbf{X}$ with i.i.d. entries from the standard Gaussian distribution and then extract its symmetric components $\mathbf{S} = \frac{1}{2}(\mathbf{X} + \mathbf{X}^T)$. Given such a random symmetric $\mathbf{S} \in \mathbb{R}^{n \times n}$ we can always approximate it by $\mathbf{\bar{S}}$ created via Algorithm 1, using $h$ Householder reflectors as in \eqref{eq:thesbar} such that
	\begin{equation}
	\mathbb{E}[\| \mathbf{S} - \mathbf{\bar{S}} \|_F^2] \leq \sum_{i=h+1}^n \sigma_i^2 - \frac{n-h}{2}.
	\label{eq:result4}
	\end{equation}
	
	\noindent \textit{Proof.} We use the eigenvalue decomposition $\mathbf{S} = \mathbf{U} \text{diag}(\mathbf{s}) \mathbf{U}^T$ and the singular value decomposition $\mathbf{S} = \mathbf{U}  \text{diag}(\mathbf{\sigma}) \mathbf{V}^T$ with the vector $\mathbf{\sigma} \in \mathbb{R}_+^n$ that has the singular values in decreasing order. We have that $\mathbf{V}^T = \mathbf{D}_\mathbf{\Sigma} \mathbf{P} \mathbf{U}^T$ where $\mathbf{D}_\mathbf{\Sigma}  = \text{diag}(\text{sign}(\mathbf{s}))$ because the singular values of $\mathbf{S}$ are the absolute values of the eigenvalues ($\sigma_i = \text{sign}(s_i) s_i$) and $\mathbf{P}$ is a permutation matrix because while the ordering of the eigenvalues in $\mathbf{s}$ is unimportant, the singular values are sorted in decreasing order, i.e., $\sigma_1 \geq \dots \geq \sigma_n$.
	
	We consider the Householder reflectors $\mathbf{J}_k$ that start the diagonalization process for $\mathbf{S}$, but only for the first $h$ steps of the process \cite[Chapter~8]{Golub1996}, i.e., a partial or incomplete eigenvalue decomposition applied to $\mathbf{S}$, as
	\begin{equation}
	\mathbf{J}_h \dots \mathbf{J}_1 \mathbf{S}  \mathbf{J}_1 \dots \mathbf{J}_h = \begin{bmatrix}
	\mathbf{\Lambda} & \mathbf{0}_{h \times (n-h)} \\
	\mathbf{0}_{(n-h) \times h} & \mathbf{\tilde{S}}
	\end{bmatrix},
	\label{eq:PartialEigenvalueStructure}
	\end{equation}
	with $\mathbf{\Lambda} = \text{diag}(\mathbf{\lambda}) \in \mathbb{R}^{h \times h}$ is a diagonal matrix with elements $\lambda_i$ such that $|\lambda_i| = \sigma_i,\ i = 1,\dots,h,$ and $\mathbf{\tilde{S}} \in \mathbb{R}^{(n-h) \times (n -h)}$ is a random symmetric matrix whose singular values are $\sigma_{h+1},\dots,\sigma_n$ (i.e, the lowest $n-h$ singular values of $\mathbf{S}$ and therefore $\| \mathbf{\tilde{S}} \|_F^2 = \sum_{i=h+1}^n \sigma_i^2 $). With this choice of the reflectors, our objective function becomes
	\begin{equation}
	\begin{aligned}
	\| \mathbf{S} -   \mathbf{J}_1 \dots \mathbf{J}_h & \text{diag}(\mathbf{\bar{s}}) \mathbf{J}_h \dots \mathbf{J}_1 \|_F^2 \\
	= & \| \mathbf{J}_h \dots \mathbf{J}_1 \mathbf{S}  \mathbf{J}_1 \dots \mathbf{J}_h - \text{diag}(\mathbf{\bar{s}}) \|_F^2,
	\end{aligned}
	\end{equation}
	that we minimize by choosing the spectra of the approximation to be $\mathbf{\bar{s}} = \begin{bmatrix}
	\mathbf{\lambda} \\ \text{diag}(\mathbf{\tilde{S}})
	\end{bmatrix}$. With this choice we have that
	\begin{equation}
	\begin{aligned}
	\mathbb{E}[\| & \mathbf{S} - \mathbf{\tilde{S}} \|_F^2] \! = \! \mathbb{E} \!\! \left[ \! \left\| \! \begin{bmatrix}
	\mathbf{0}_{h \times h} & \mathbf{0}_{h \times (n-h)} \\
	\mathbf{0}_{(n-h) \times h} & \mathbf{\tilde{S}} \! - \! \text{diag}(\text{diag}(\mathbf{\tilde{S}}))
	\end{bmatrix} \! \right\|_F^2 \! \right] \\
	= & \mathbb{E}[\|\mathbf{\tilde{S}} - \text{diag}(\text{diag}(\mathbf{\tilde{S}}))  \|_F^2] \\
	= & \mathbb{E}[\| \mathbf{\tilde{S}}  \|_F^2 + \| \text{diag}(\mathbf{\tilde{S}}) \|_2^2 - 2\text{tr}(\mathbf{\tilde{S}}^T \text{diag}(\text{diag}(\mathbf{\tilde{S}})))] \\
	= & \sum_{i=h+1}^n \sigma_i^2 - \frac{n-h}{2} \approx \frac{(n-h)^2 - (n-h)}{2}.
	\end{aligned}
	\label{eq:bigformula}
	\end{equation}
	For the final equality, which is accurate for large $n$ and $h$, we used $\mathbb{E}[\tilde{S}_{ij}^2] = \frac{1}{\sqrt{2}}$ and that we are computing the Frobenius norm of a matrix with $(n-h)(n-h-1)$ non-zero elements.
	
	This is an upper bound to the approximation accuracy of Algorithm 1 since we could initialize the procedure with the $h$ reflectors that achieve \eqref{eq:PartialEigenvalueStructure} and then iteratively improve them by monotonically decreasing the objective function. The result also shows that the proposed approximation is at least as good as the low-rank approximation of $\mathbf{S}$.$\hfill \blacksquare$
	
	\noindent\textbf{Remark (the case of positive semidefinite $\mathbf{S}$).} Given a random matrix $\mathbf{X}$ with entries i.i.d. from the standard Gaussian distribution we call the symmetric positive definite $\mathbf{S} = \mathbf{X}\mathbf{X}^T$ a Wishart matrix which is diagonally dominant (we have $\mathbb{E}[S^2_{ii}] = n^2 + 2n$ and $\mathbb{E}[S^2_{ij}] = n,\ i \neq j$). When $n$ is large and $h \ll n$ we also expect that $\mathbf{\tilde{S}}$ in \eqref{eq:bigformula} is diagonally dominant and therefore it is well approximated by $\text{diag}(\text{diag}(\mathbf{\tilde{S}}))$ leading to a lower approximation error as compared to the indefinite case when $\mathbf{S} = \frac{1}{2}(\mathbf{X} + \mathbf{X}^T)$. $\hfill \blacksquare$

	\section{Results}
	
	In this section we show the approximation performance of the two proposed reflector structures, for orthonormal and symmetric matrices. To measure the accuracy of our approximations, given a target matrix $\mathbf{X}$ and its approximation $\mathbf{\bar{X}}$ we use the normalized relative representation error:
	\begin{equation}
	\epsilon(\mathbf{X}, \mathbf{\bar{X}}) = \frac{1}{4}\frac{\| \mathbf{X} - \mathbf{\bar{X}} \|_F^2}{\| \mathbf{X} \|_F^2}.
	\label{eq:reperror}
	\end{equation}
	This error is normalized by $4 \| \mathbf{X} \|_F^2$ such that for any pair of orthonormal matrices $(\mathbf{X}, \mathbf{\bar{X}})$ we have that $0 \leq \epsilon(\mathbf{X}, \mathbf{\bar{X}}) \leq 1$. 
	
	\subsection{Synthetic experiments: orthonormal case}
	
	In this section, we randomly generate orthonormal $\mathbf{U}$: we construct a random matrix with entries from the standard Gaussian distribution on which we apply the QR decomposition and keep the orthonormal component.

\begin{figure*}[!tbp]
		\centering
		\begin{minipage}[t]{0.32\textwidth}
			\includegraphics[trim = 18 5 30 15, clip, width=\textwidth]{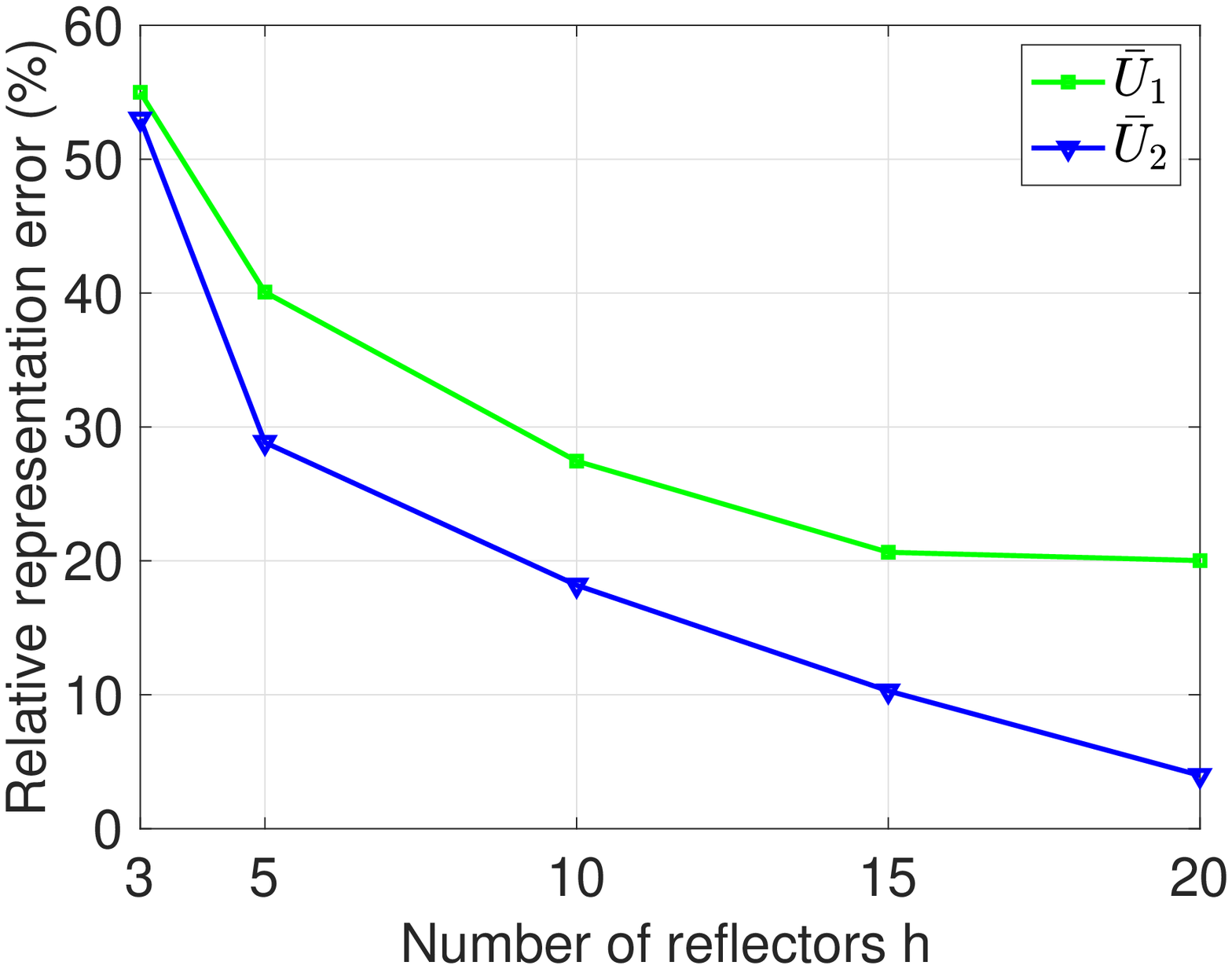}
			\caption{For fixed $n = 32$, we show the relative representation error \eqref{eq:reperror} as a function of the number of reflectors $h$ for the two approximations we consider: $\mathbf{\bar{U}}_1$ and $\mathbf{\bar{U}}_2$ as in \eqref{eq:approxWithReflectors} and \eqref{eq:approxWithReflectors2}, respectively. The results are averaged over 100 realizations of random orthonormal matrices in each case. Complexity of unstructured matrix-vector multiplication is reached for $h=16$.}
			\label{fig:householder0}
		\end{minipage}
		\hfill
		\begin{minipage}[t]{0.32\textwidth}
			\includegraphics[trim = 18 5 20 15, clip, width=\textwidth]{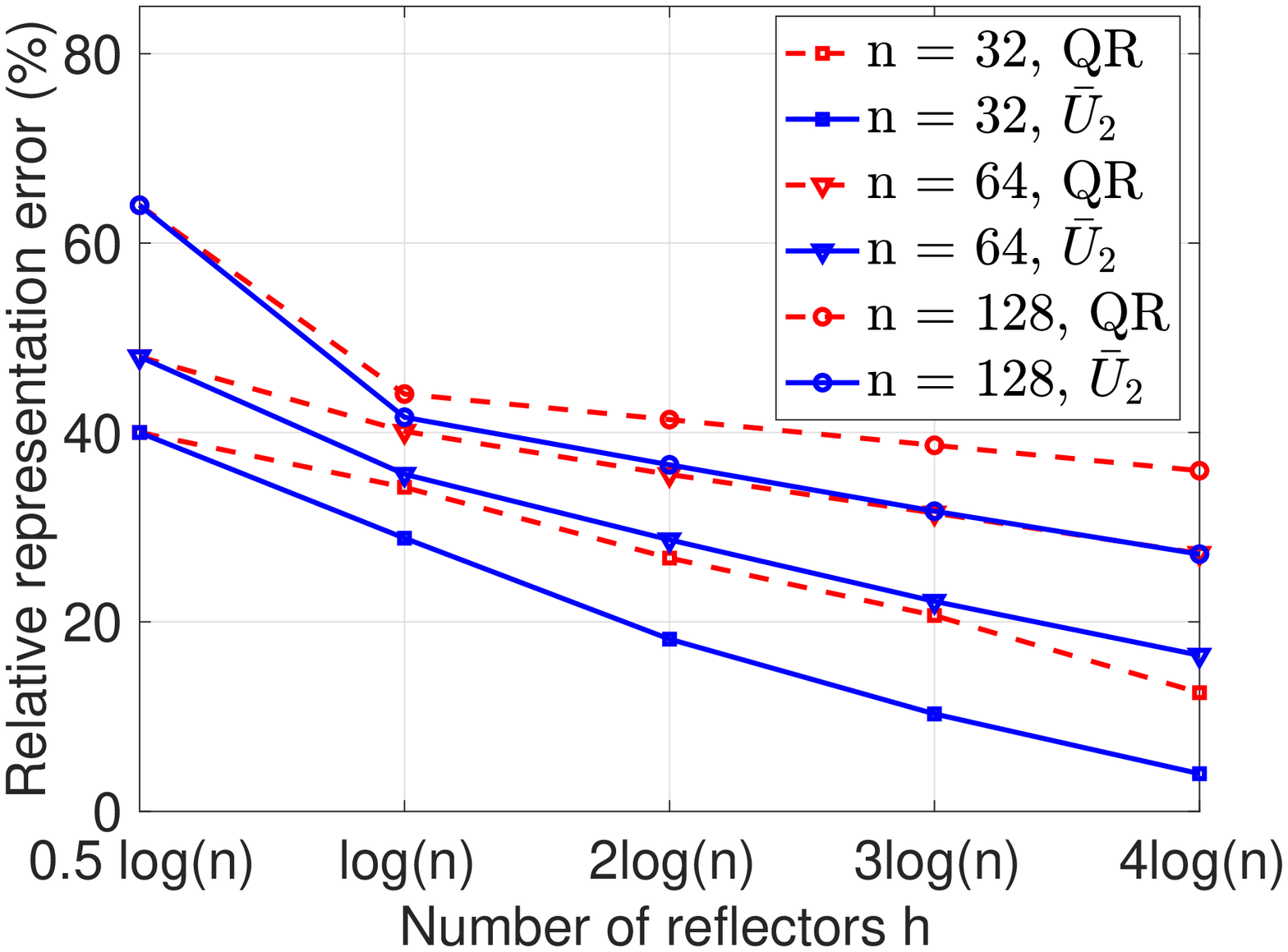}
			\caption{Relative representation error \eqref{eq:reperror} achieved by the proposed approximation $\mathbf{\bar{U}}_2$ as a function of the number of reflectors $h$ for varying dimensions $n \in \{32, 64, 128\}$. The dotted red lines with the corresponding symbols show the same type of results achieved by partial QR decomposition \eqref{eq:ReflectorStructure}.}
			\label{fig:householder1}
		\end{minipage}
		\hfill
		\begin{minipage}[t]{0.32\textwidth}
			\includegraphics[trim = 18 5 30 15, clip, width=\textwidth]{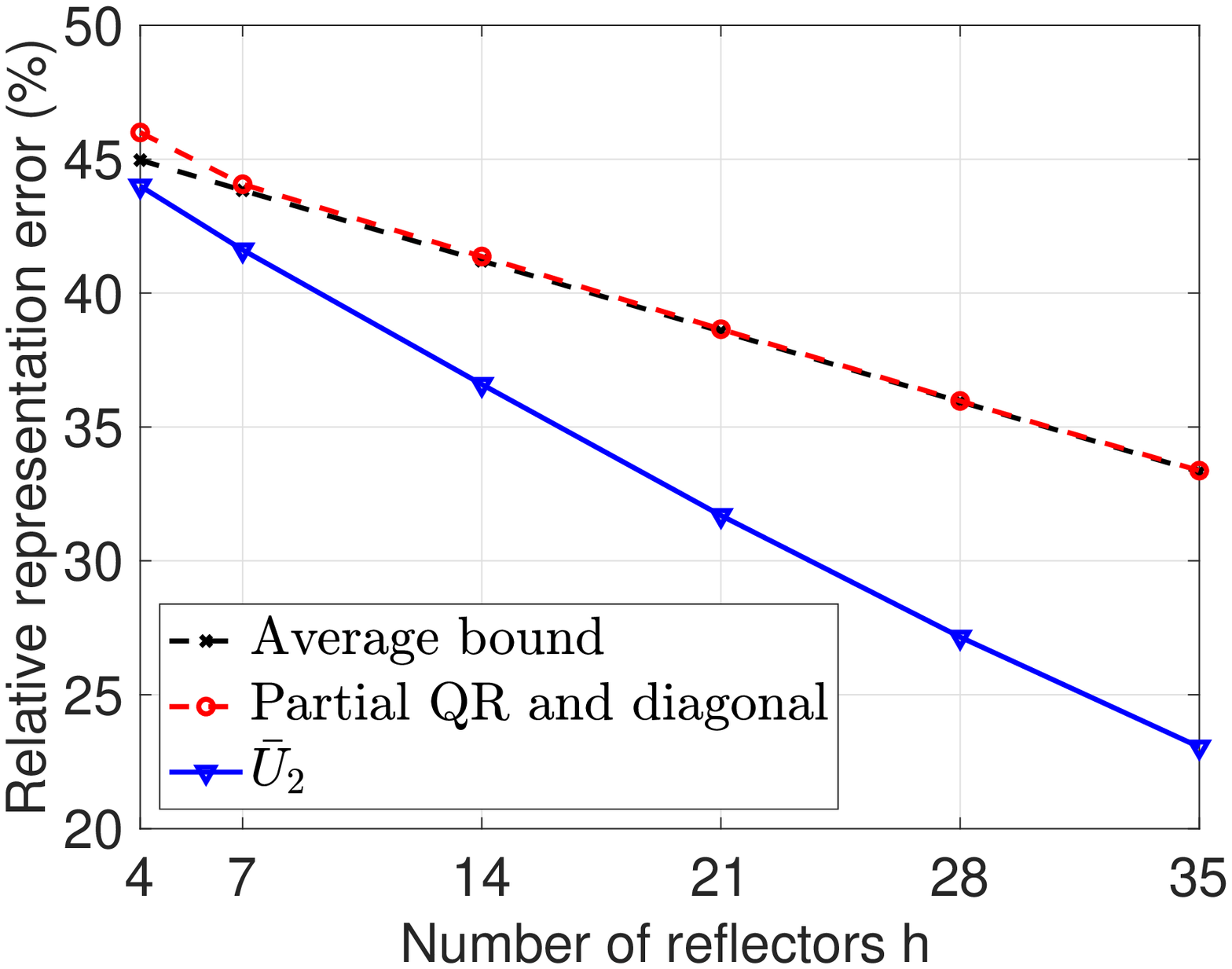}
			\caption{Relative representation error \eqref{eq:reperror} achieved by the proposed method as a function of the number of reflectors $h$ and dimension $n=128$. The dashed red line shows the same type of results achieved by partial QR and diagonal decomposition in \eqref{eq:ReflectorStructure} and the black dashed line (almost invisible due to the overlap with the previous line) shows the average bound of Result 3.}
			\label{fig:householder2}
		\end{minipage}
	\end{figure*}

	In Figure \ref{fig:householder0} we compare the two orthonormal approximations we consider in this paper, $\mathbf{\bar{U}}_1$ and $\mathbf{\bar{U}}_2$. As discussed, even from a theoretical perspective, the approximation structure $\mathbf{\bar{U}}_2$ outperforms (or matches, in the worst case) $\mathbf{\bar{U}}_1$ since the former has more degrees of freedom in its model. We also confirmed this fact by the simulation results in Figure \ref{fig:householder0} where the gap between the two is approximately $10\%$ in the relative representation error. Furthermore, given enough reflectors $h$, at most $n-1$ constructed by the QR decomposition, the approximation $\mathbf{\bar{U}}_2$ can reach zero representation error while $\mathbf{\bar{U}}_1$ can rarely achieve perfect reconstruction regardless of the number of reflectors we allow.
	
	In Figures \ref{fig:householder1} and \ref{fig:householder2} we show the relative representation error obtained when approximating a random orthonormal $\mathbf{U}$ by $\mathbf{\bar{U}}_2$ as in \eqref{eq:approxWithReflectors2}. In Figure \ref{fig:householder2} we also show the average bound developed in Result 1, which serves here as an upper bound for the proposed approach. As the number of reflectors $h$ increases, the bound is essentially tight and matches the approximation accuracy of the $h$-step QR factorization. Also, notice that the difference between our proposed method and the bound increases with the number of reflectors $h$.
	
	The approximations $\mathbf{\bar{U}}_1$ and $\mathbf{\bar{U}}_2$ apply also to classic transforms. If $\mathbf{F} \in \mathbb{C}^{n \times n}$ is the Fourier matrix then $\mathbf{F} + \mathbf{F}^H$ has $\left\lfloor \frac{n+2}{4} \right\rfloor$ eigenvalues with value $-2$ and $\left\lfloor \frac{n+4}{4} \right\rfloor$ eigenvalues with value $2$ and the others up to $n$ are zero. For the Hadamard matrix $\mathbf{H} \in \mathbb{R}^{n \times n}$ we have that $\mathbf{H} + \mathbf{H}^T$ has half the eigenvalues with value $-2$ and half with value $2$. For both these transforms, the approximation accuracies of $\mathbf{\bar{U}}_1$ and $\mathbf{\bar{U}}_2$ are identical. For other well known transforms, like the discrete cosine matrix (which is related to the Fourier) and the Haar matrix, the eigenvalues $z_k$ do not all have extremal values $\{ \pm 2 \}$ and therefore the approximation analyses with the proposed $\mathbf{\bar{U}}_1$ and $\mathbf{\bar{U}}_2$ cannot be easily done analytically. Regardless, we note that the proposed factorizations are not appropriate to build numerically efficient and accurate approximations of these structures (according to Results 1 and 2 that highlight the need of a small number of extreme eigenvalues in the matrices to be approximated). 

	\subsection{Synthetic experiments: symmetric case}
	
	In this section, we randomly generate symmetric $\mathbf{S}$: we construct a random matrix $\mathbf{X}$ with entries from the standard Gaussian distribution and compute either $\mathbf{X}+\mathbf{X}^T$ or $\mathbf{XX}^T$.
	
\begin{figure*}[!tbp]
		\centering
		\begin{minipage}[t]{0.32\textwidth}
			\includegraphics[trim = 18 0 20 15, clip, width=\textwidth]{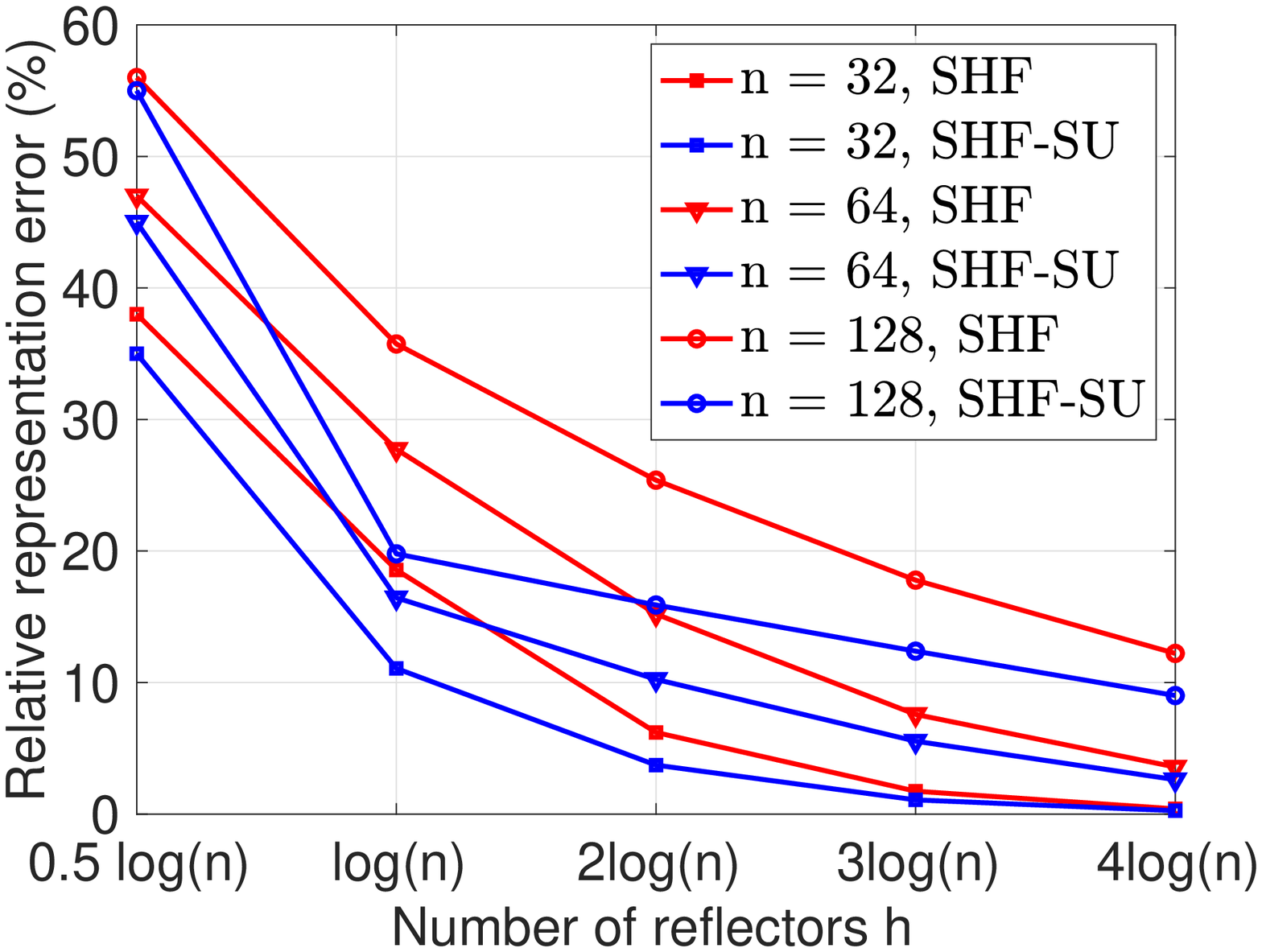}
			\caption{Relative representation error \eqref{eq:reperror} achieved by the proposed approximation $\mathbf{\bar{S}}$ as a function of the number of reflectors $h$ for varying dimensions $n \in \{32, 64, 128\}$. The red lines with the corresponding symbols show the same type of results achieved by the same algorithm where we also allow the spectrum update via \eqref{eq:updateofsbar}. The results are averaged over 100 realizations of random symmetric matrices $\mathbf{S} = \frac{1}{2}(\mathbf{X} + \mathbf{X}^T)$ in each case, where $\mathbf{X}$ is a matrix with entries i.i.d. standard Gaussian.}
			\label{fig:householdersym1}
		\end{minipage}
		\hfill
		\begin{minipage}[t]{0.32\textwidth}
			\includegraphics[trim = 18 0 20 15, clip, width=\textwidth]{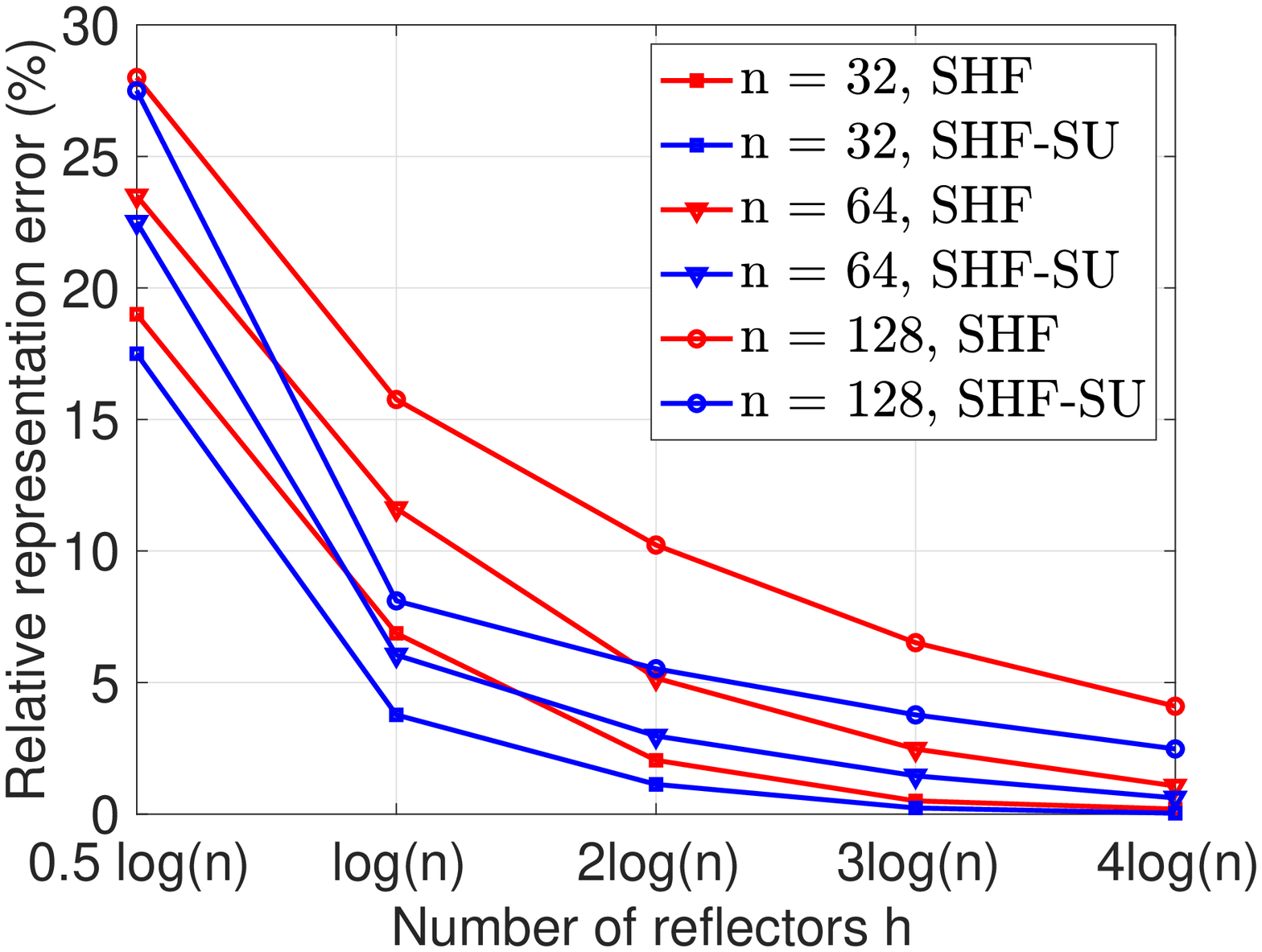}
			\caption{Experimental results analogous to Figure \ref{fig:householdersym1} for a symmetric positive definite $\mathbf{S} = \mathbf{XX}^T$. In both figures, complexity of unstructured matrix-vector multiplication is reached when: $h = \lfloor 3.2\log n\rfloor$ for $n = 32$, $h  = \lfloor5.3 \log n\rfloor$ for $n = 64$ and $h=\lfloor 9.1\log n \rfloor$ for $n =  128$.}
			\label{fig:householdersym2}
		\end{minipage}
		\hfill
		\begin{minipage}[t]{0.32\textwidth}
			\includegraphics[trim = 18 0 22 15, clip, width=\textwidth]{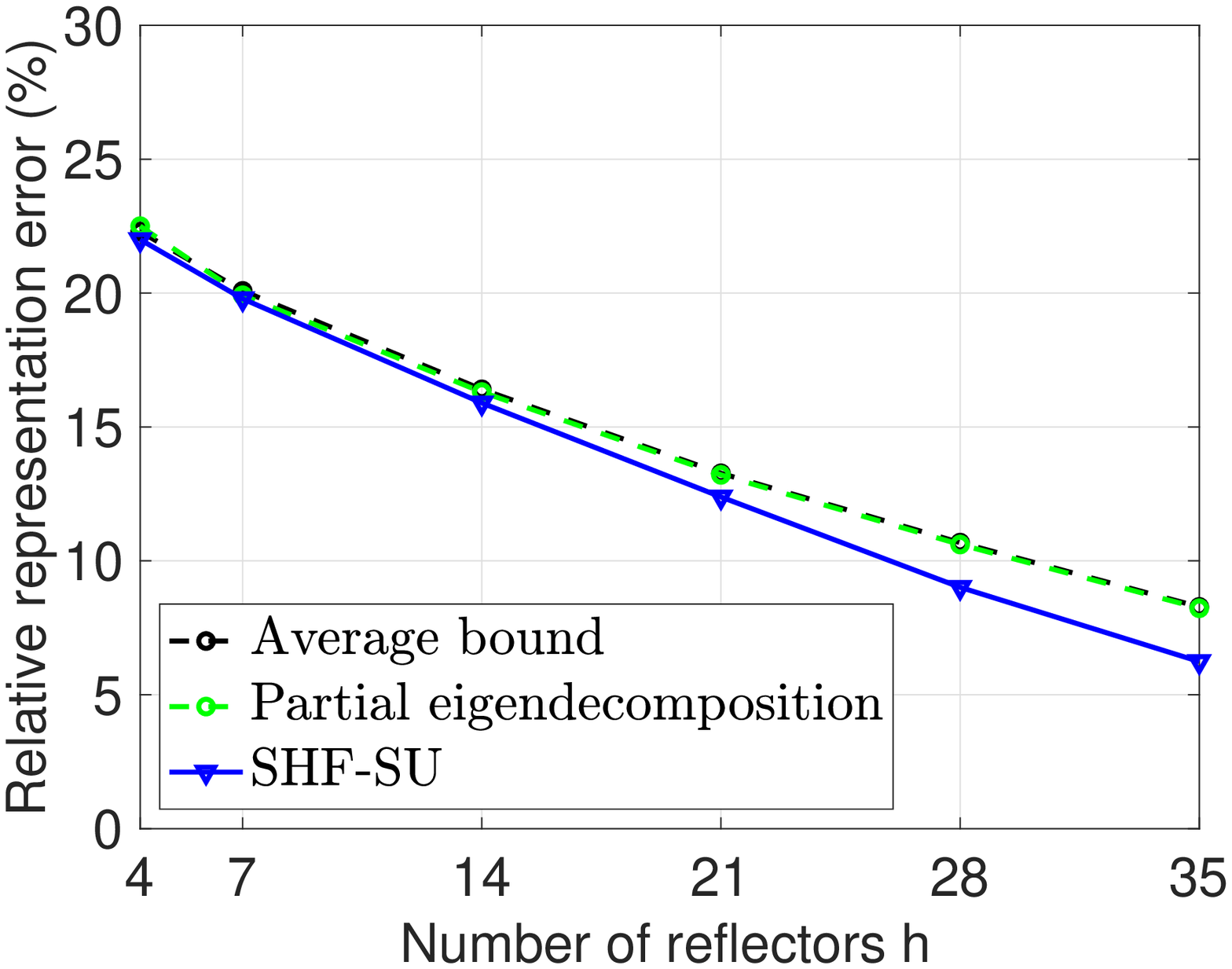}
			\caption{Relative representation error \eqref{eq:reperror} achieved by SHF--SU as a function of the number of reflectors $h$ and dimension $n=128$. The dashed green line shows the same type of results achieved by partial eigendecomposition and diagonal update in \eqref{eq:PartialEigenvalueStructure} and the black dashed line (almost invisible due to the overlap with the previous line) shows the average bound of Result 4.}
			\label{fig:householdersym4}
		\end{minipage}
	\end{figure*}

	Figures \ref{fig:householdersym1} and \ref{fig:householdersym2} show experimental simulations for approximating random symmetric matrices indefinite and positive definite, respectively. The first observation is that, as expected, the approximation accuracy with the proposed structure is better when considering positive definite symmetric matrices (Figure \ref{fig:householdersym2} shows lower representation errors as compared to Figure \ref{fig:householdersym1} for the same dimension $n$). The second observation is that, in both situations, allowing for spectrum updates significantly increases the approximation accuracy: on average SHF--SU achieves approximately twice the representation accuracy as compared to SHF.
	
\begin{figure*}[!tbp]
		\centering
		\begin{minipage}[t]{0.32\textwidth}
			\includegraphics[trim = 5 0 30 15, clip, width=\textwidth]{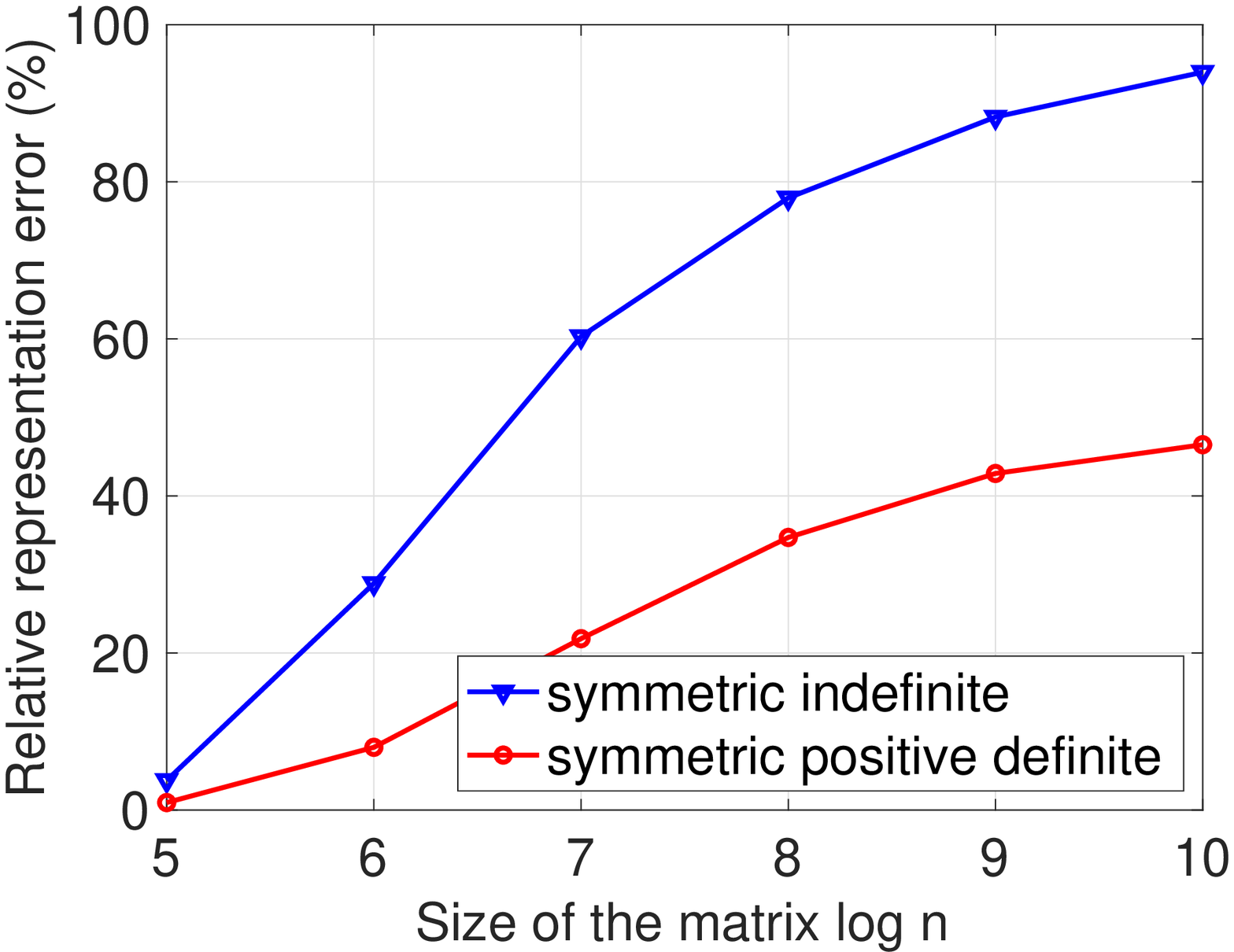}
			\caption{Relative representation error \eqref{eq:reperror} achieved by SHF--SU for various $n$ and fixed $h = 16$ for both indefinite and positive definite matrices.}
			\label{fig:final_results}
		\end{minipage}
		\hfill
		\begin{minipage}[t]{0.32\textwidth}
			\includegraphics[trim = 18 5 30 15, clip, width=\textwidth]{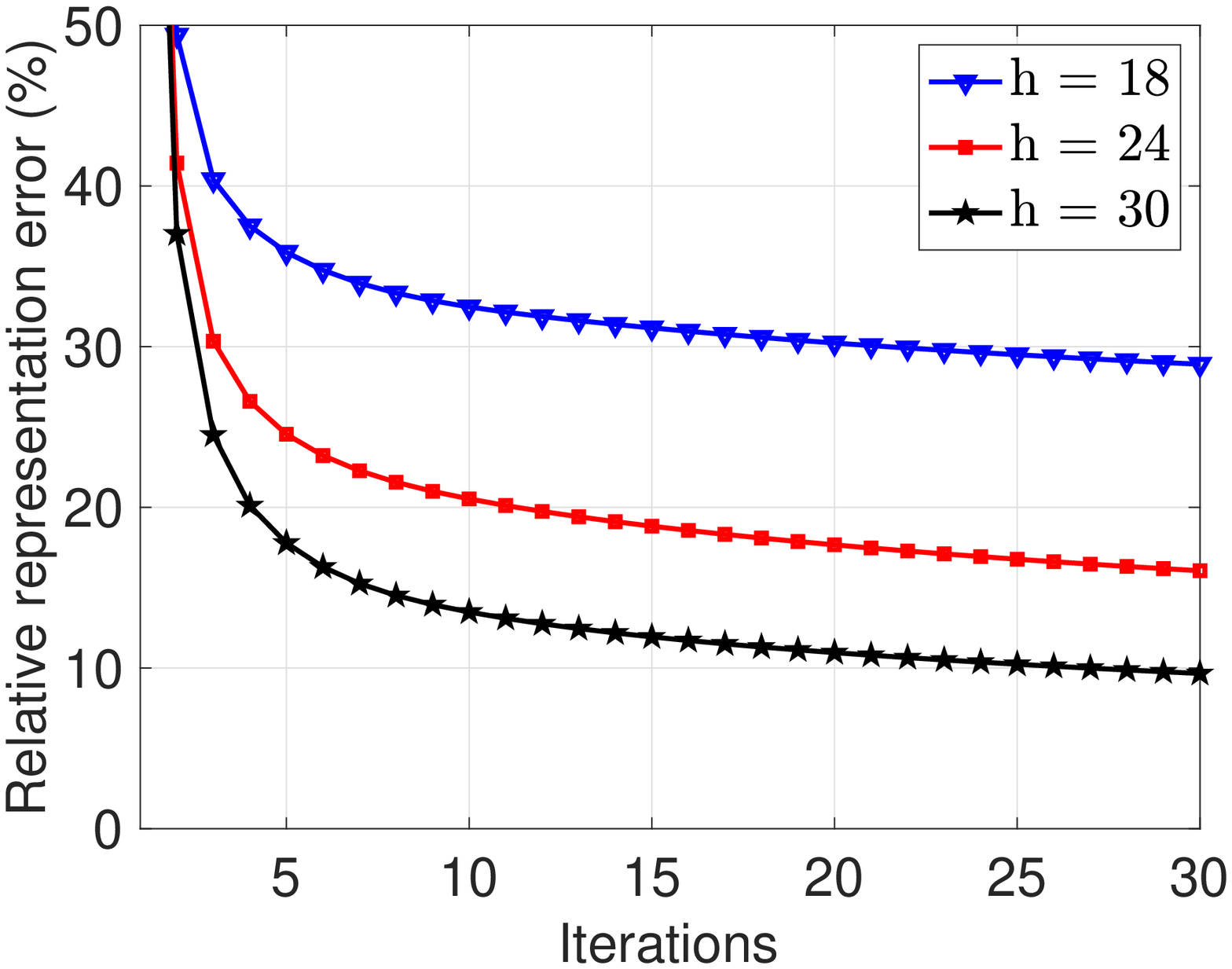}
			\caption{For fixed $n = 64$ and reflectors $h \in \{ 18, 24, 30 \}$, we show the evolution of the relative representation error as a function of the number of iterations in SHF. The results are averaged over 100 realizations random indefinite matrices.}
			\label{fig:plot_iter_symm}
		\end{minipage}
		\hfill
		\begin{minipage}[t]{0.32\textwidth}
			\includegraphics[trim = 18 5 20 15, clip, width=\textwidth]{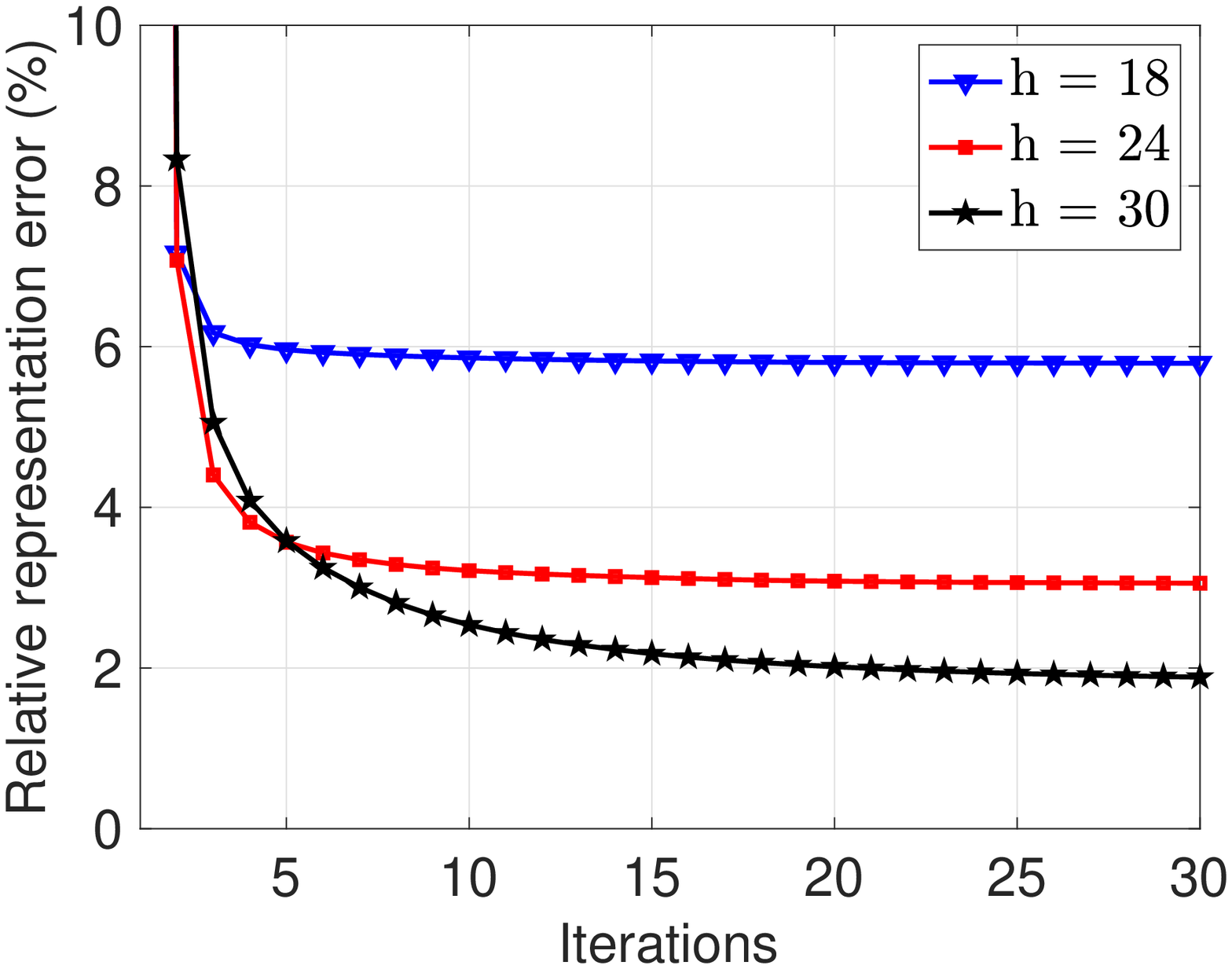}
			\caption{Experimental results analogous to Figure \ref{fig:plot_iter_symm} for symmetric positive definite matrices.}
			\label{fig:plot_iter_posdef}
		\end{minipage}
	\end{figure*}

	In Figure \ref{fig:householdersym4} we show the approximation gap between the bound developed in Result 4 and the proposed SHF--SU algorithm for $n=128$. As expected, SHF--SU always performs better than the eigendecomposition approach. The gap is not large due to the eigenvalue distribution of the randomly generated $\mathbf{S}$, i.e., most of the energy is concentrated in a few (highest) eigenvalues that are also captured well in the $h$ eigenvalue decomposition. In Figure \ref{fig:final_results} we apply SHF--SU with fixed number of reflectors $h=16$ and we increase $n$ to observe the deterioration of the representation error while we separate again the indefinite and positive definite cases. Figures \ref{fig:plot_iter_symm} and \ref{fig:plot_iter_posdef} show the evolution of the SHF--SU with each iteration for the symmetric indefinite case and for the positive definite case, respectively. For ease of exposition, we show only the first 30 iterations in both cases but notice that while for the positive definite case we already get convergence (progress is below $10e-5$) in the indefinite case, even after 30 iterations there is still significant progress. We observe this behavior in all cases, with indefinite matrices SHF or SHF--SU converges in more iterations (more than double, on average).
	
	We would like to highlight that the choice of the number of Householder reflectors $h$ in the approximation is not obvious. In the ideal situation, we would know what computational budget is available in the application as this would lead to a possible estimation of $h$. Alternatively, a Pareto curve showing the trade-off between computational complexity and approximation error (relative representation error versus $h$) could help the human operator reach a decision. Finally, choosing $h$ might also be done using information-theoretic criteria for model selection.
	
	\subsection{Application: learning fast distance metric transformations}
	
	In the context of machine learning algorithms, the classification accuracy of many methods significantly depends on the choice of a good metric, i.e., a good distance between any two data points of the dataset. For example, the performance of the $k$-nearest neighbors algorithm ($k$-NN) \cite{kNN} highly depends on using a metric that accurately reflects the relationship between data points (both, data points from the same class as well as data points from different classes). The simple, standard Euclidean distance regularly used by $k$-NN does not exploit any possible structure that might exist in the data. Given labeled data points a problem that arises is constructing a distance metric such that points from the same class are ``close" and points from different classes are ``far". This is known as the distance metric learning problem \cite{SurveyMetricLearning}: given a training dataset, find a linear transformation of the input data such that points from the same class are concentrated while the separation between points of different classes increases. This technique has been shown to consistently produce improved results as compared to the Euclidean distance \cite{Shental2002, ShalevShwartz2004, Weinberger2005}.
	
	Concretely, consider a labeled training set $\{\mathbf{x}_i, y_i\}_{i=1}^N$ with inputs $\mathbf{x}_i \in \mathbb{R}^n$ and discrete, finite class labels $\mathbf{y}_i$. Instead of using the Euclidean distance between points, i.e., $d(\mathbf{x}_i, \mathbf{x}_j) = \| \mathbf{x}_i - \mathbf{x}_j \|_2^2 = (\mathbf{x}_i - \mathbf{x}_j)^T (\mathbf{x}_i - \mathbf{x}_j)$, our goal is to learn a linear transformation $\mathbf{L} \in \mathbb{R}^{n \times n}$ such that we use the new distance $d_\mathbf{S} (\mathbf{x}_i, \mathbf{x}_j) = \| \mathbf{L}(\mathbf{x}_i - \mathbf{x}_j) \|_2^2= (\mathbf{x}_i - \mathbf{x}_j)^T \mathbf{L}^T \mathbf{L} (\mathbf{x}_i - \mathbf{x}_j)$. If we denote $\mathbf{S} = \mathbf{L}^T \mathbf{L}$, then the problem of learning the symmetric positive semidefinite metric $\mathbf{S} \in \mathbb{R}^{n \times n}$ has the name of Mahalanobis metric learning. The metric $\mathbf{S}$ is optimized with the goal that nearest neighbors always belong to the same class while examples from different classes are maximally separated. We can write the following optimization problem:
	\begin{equation}
	\underset{\mathbf{S} \in \mathbb{S}_+^{n \times n}}{\text{maximize}} \! \! \! \sum_{ (i, j) \in \mathcal{D} } \! \! \! \!  d_\mathbf{S} (\mathbf{x}_i, \mathbf{x}_j)
	\text{ subject to} \! \! \! \sum_{ (i, j) \in \mathcal{S} } \! \! \! \! d_\mathbf{S} (\mathbf{x}_i, \mathbf{x}_j) \leq 1,
	\end{equation}
	where $\mathbb{S}_+^{n \times n}$ is the set of $n \times n$ symmetric positive semidefinite matrices, $\mathcal{D}$ and $\mathcal{S}$ are sets of the pairs of points that belong to different or similar classes, respectively.

	Given the learned $\mathbf{S}$ and a new data point $(\mathbf{x}_{N+1}, y_{N+1})$ the first step will be to perform the transformation $\mathbf{x}_{N+1} \leftarrow \mathbf{L} \mathbf{x}_{N+1}$, before running the $k$-NN algorithm. The transform is taken to be $\mathbf{L} = \sqrt{\mathbf{\Lambda}} \mathbf{U}$ where $\mathbf{S} = \mathbf{U}^T \mathbf{\Lambda U}$ is the eigendecomposition of $\mathbf{S}$. In general, the cost of performing $\mathbf{L} \mathbf{x}_{N+1}$ is $O(n^2)$ and our goal is to reduce it to $O(n \log n)$ by using the Householder reflector factorization proposed in this paper for the orthonormal matrix $\mathbf{U}$.
	
	In this section, we follow the distance metric learning for large margin $k$-NN classification algorithm (LMNN) \cite{Weinberger2005}. Given a total budget of $K = 100$ iterations, we run the regular algorithm for 50 iterations and then for the latter we factor the learned distance metric using SHF--SU with a given $h$, at each iteration. The idea is to build an unconstrained distance metric and then project it on the set of fast transformations built with $h$ Householder reflectors. We use and modify the LMNN metric learning implementation provided in the Matlab Toolbox for Dimensionality Reduction\footnote{http://homepage.tudelft.nl/19j49}. We test the approach using two well-known test datasets: ISOLET\footnote{https://archive.ics.uci.edu/ml/datasets/isolet} with $7797$, NEWS\footnote{http://people.csail.mit.edu/jrennie/20Newsgroups} with $18828$ and MNIST\footnote{http://yann.lecun.com/exdb/mnist/} with $70000$ data points, respectively. In both cases, we use principal component analysis (PCA) to reduce the dimensionality of the data, to $n=100$, $n=200$ and $n=40$ components, respectively. Then we use a 70/30 split of the data for training and testing and we average the accuracy of the classification over 100 random realizations. Results are shown in Table \ref{table:results}.
	
	Notice that the average accuracy is very close to that of the general (unconstrained) metric while the cost of $\mathbf{Lx}_{N+1}$ is now controllable, i.e., $4nh$ operations instead of $2n^2$. For example, when $h = \lceil \log_2 n \rceil$ the speedup at test time due to the proposed factorization is approximately $\times 7$, $\times 4$, $\times 2$, $\times 12$ and $\times 3$, respectively for each dataset as ordered in Table \ref{table:results}. Notice that at training time we have an extra computational overhead added by constructing the factorization. This overhead is not significant in general, as it depends only on the dimension $n$ and not the size of the overall training dataset, $N$.
	\begin{table}[t]
		\begin{center}
			\caption{Average testing error of numerically efficient large margin $k$-NN metric learning with Householder reflectors.}\label{table:results}
			\begin{tabular}{|c||c|c|c|c|}
				\hline
				\multicolumn{2}{|c}{} & \multicolumn{3}{|c|}{Number of reflectors $h$} \\ \hline
				Dataset & Full metric \cite{Weinberger2005} &  $\lceil \log_2n \rceil$ & $ \lceil 2\log_2 n \rceil$ & $\lceil 3\log_2 n \rceil$ \\ \hline \hline
				ISOLET & \textbf{4.5\%} & 8.7\% & 6.4\% & 4.6\% \\ \hline
				USPS & \textbf{4.0\%} & 7.5\% & 5.6\% & 4.0\% \\ \hline
				UCI & \textbf{8.0\%} & 13.1\% & 9.4\% & 8.2\% \\ \hline
				NEWS & \textbf{13.1\%} & 18.2\% & 14.1\% & 13.3\% \\ \hline
				MNIST & \textbf{1.9\%} & 4.5\% & 2.3\% & 2.0\% \\ \hline
			\end{tabular}
		\end{center}
	\end{table}

	\section{Conclusions}
	
	In this paper, we describe a class of orthonormal matrices that are the product of just a few Householder reflectors. By controlling the number of reflectors in the factorization, we regulated the computational complexity of matrix-vector multiplications with these transformations. We perform an analysis of the proposed structures and describe algorithms that approximate (imperfectly in general) any orthonormal operator by a product of a given number of Householder reflectors. We then propose a similar factorization for symmetric matrices and we show an application in the context of a $k$-nearest neighbors classification problem where we use the proposed factorizations to approximate a learned distance metric with little performance degradation in terms of the classification accuracy and little computational overhead in the training phase but with a significant computational speedup for the testing phase.

	\bibliographystyle{IEEEtran}
	\bibliography{refs}

\begin{thebibliography}{10}
\providecommand{\url}[1]{#1}
\csname url@samestyle\endcsname
\providecommand{\newblock}{\relax}
\providecommand{\bibinfo}[2]{#2}
\providecommand{\BIBentrySTDinterwordspacing}{\spaceskip=0pt\relax}
\providecommand{\BIBentryALTinterwordstretchfactor}{4}
\providecommand{\BIBentryALTinterwordspacing}{\spaceskip=\fontdimen2\font plus
\BIBentryALTinterwordstretchfactor\fontdimen3\font minus
  \fontdimen4\font\relax}
\providecommand{\BIBforeignlanguage}[2]{{%
\expandafter\ifx\csname l@#1\endcsname\relax
\typeout{** WARNING: IEEEtran.bst: No hyphenation pattern has been}%
\typeout{** loaded for the language `#1'. Using the pattern for}%
\typeout{** the default language instead.}%
\else
\language=\csname l@#1\endcsname
\fi
#2}}
\providecommand{\BIBdecl}{\relax}
\BIBdecl

\bibitem{FFT}
J.~Cooley and J.~Tukey, ``An algorithm for the machine calculation of complex
  {Fourier} series,'' \emph{Mathematics of Computation}, vol.~19, no.~90, pp.
  297--301, 1965.

\bibitem{FFTComputational}
C.~F. van Loan, \emph{Computational Frameworks for the Fast Fourier
  Transform}.\hskip 1em plus 0.5em minus 0.4em\relax SIAM, 1992.

\bibitem{DCTComputation}
M.~Narasimha and A.~Peterson, ``On the computation of the discrete cosine
  transform,'' \emph{IEEE Trans. Comm.}, vol.~26, no.~6, pp. 934--936, 1978.

\bibitem{DCTComputation2}
E.~Feig and S.~Winograd, ``Fast algorithms for the discrete cosine transform,''
  \emph{IEEE Trans. Sig. Process.}, vol.~40, no.~9, pp. 2174--2193, 1992.

\bibitem{HartleyComputation}
R.~N. Bracewell, ``Computing with the {Hartley} transform,'' \emph{Computers in
  Physics}, vol.~9, no.~4, pp. 373--379, 1995.

\bibitem{WalshComputation}
B.~J. Fino and V.~R. Algazi, ``Unified matrix treatment of the fast
  {Walsh-Hadamard} transform,'' \emph{IEEE Trans. Computers}, vol.~25, no.~11,
  pp. 1142--1146, 1976.

\bibitem{FWT}
G.~Beylkin, R.~Coifman, and V.~Rokhlin, ``Fast wavelet transforms and numerical
  algorithms,'' \emph{Comm. Pure Appl. Math.}, vol.~44, pp. 141--183, 1991.

\bibitem{doi:10.1002/cpa.3160410705}
I.~Daubechies, ``Orthonormal bases of compactly supported wavelets,''
  \emph{Communications on Pure and Applied Mathematics}, vol.~41, no.~7, pp.
  909--996, 1988.

\bibitem{10.1007/978-3-642-97177-8_2}
Y.~Meyer, ``Orthonormal wavelets,'' in \emph{Wavelets}, J.-M. Combes,
  A.~Grossmann, and P.~Tchamitchian, Eds.\hskip 1em plus 0.5em minus
  0.4em\relax Berlin, Heidelberg: Springer Berlin Heidelberg, 1989, pp. 21--37.

\bibitem{HaarCompute}
A.~Haar, ``Zur {Theorie} der orthogonalen {Funktionensysteme},''
  \emph{Mathematische Annalen}, vol.~69, no.~3, pp. 331--371, 1910.

\bibitem{LiftingScheme}
W.~Sweldens, ``The lifting scheme: A construction of second generation
  wavelets,'' \emph{{SIAM} Journal on Mathematical Analysis}, vol.~29, no.~2,
  pp. 511--546, 1997.

\bibitem{Householder:1958:UTN:320941.320947}
A.~S. Householder, ``Unitary triangularization of a nonsymmetric matrix,''
  \emph{J. ACM}, vol.~5, no.~4, pp. 339--342, 1958.

\bibitem{Golub1996}
G.~H. Golub and C.~F. van Loan, \emph{Matrix Computations}.\hskip 1em plus
  0.5em minus 0.4em\relax Johns Hopkins University Press, 1996.

\bibitem{Steinhardt9259}
A.~O. Steinhardt, ``Householder transforms in signal processing,'' \emph{IEEE
  ASSP Magazine}, vol.~5, no.~3, pp. 4--12, 1988.

\bibitem{DictHouseholder}
C.~Rusu, N.~Gonzalez-Prelcic, and R.~Heath, ``Fast orthonormal sparsifying
  transforms based on {Householder} reflectors,'' \emph{IEEE Trans. Sig.
  Process.}, vol.~64, no.~24, pp. 6589--6599, 2016.

\bibitem{lemagoarou:hal-01104696}
L.~Le~Magoarou and R.~Gribonval, ``{Chasing butterflies: In search of efficient
  dictionaries},'' in \emph{ICASSP}, Brisbane, Australia, 2015.

\bibitem{FastSparsifyingTransforms}
C.~Rusu and J.~Thompson, ``Learning fast sparsifying transforms,'' \emph{IEEE
  Transactions on Signal Processing}, vol.~65, no.~16, pp. 4367--4378, 2017.

\bibitem{MatrixAnalysis}
R.~A. Horn and C.~R. Johnson, \emph{Matrix analysis}.\hskip 1em plus 0.5em
  minus 0.4em\relax Cambridge University Press, 2013.

\bibitem{CollinsMale2011}
B.~Collins and C.~Male, ``The strong asymptotic freeness of {Haar} and
  deterministic matrices,'' \emph{Annales Scientifiques de l'Ecole Normale
  Superieure}, vol.~47, 2011.

\bibitem{OrthoTheory}
E.~Borel, \emph{Introduction g\'{e}om\'{e}trique $\grave{a}$ quelques
  th\'{e}ories physiques}.\hskip 1em plus 0.5em minus 0.4em\relax
  Gauthier-Villars, Paris, 1906.

\bibitem{SparsePCA}
H.~Zou, T.~Hastie, and R.~Tibshirani, ``Sparse principal component analysis,''
  \emph{Journal of Computational and Graphical Statistics}, vol.~15, no.~2, pp.
  262--286, 2006.

\bibitem{OnTheKron}
K.~Schacke, ``On the {Kronecker} product,'' \emph{Masters Thesis, University of
  Waterloo}, 2013.

\bibitem{ubiquitousKronecker}
C.~F. Van~Loan, ``The ubiquitous {Kronecker} product,'' \emph{J. Computational
  and Applied Mathematics}, vol. 123, pp. 85--100, 2000.

\bibitem{ProdQuadForms}
M.~Lin and G.~Sinnamon, ``A condition for convexity of a product of positive
  definite quadratic forms,'' \emph{SIAM J. Matrix Anal. Appl.}, vol.~32,
  no.~2, pp. 457--462, 2011.

\bibitem{10.2307/j.ctt7smmk}
P.-A. Absil, R.~Mahony, and R.~Sepulchre, \emph{Optimization Algorithms on
  Matrix Manifolds}.\hskip 1em plus 0.5em minus 0.4em\relax Princeton
  University Press, 2008.

\bibitem{kNN}
T.~Cover and P.~Hart, ``Nearest neighbor pattern classification,'' \emph{IEEE
  Transactions in Information Theory}, vol. IT-13, pp. 21--27, 1967.

\bibitem{SurveyMetricLearning}
A.~Bellet, H.~Amaury, and S.~Marc, ``A survey on metric learning for feature
  vectors and structured data,'' \emph{arXiv:1306.6709}, 2013.

\bibitem{Shental2002}
N.~Shental, T.~Hertz, D.~Weinshall, and M.~Pavel, ``Adjustment learning and
  relevant component analysis,'' in \emph{Seventh European Conference on
  Computer Vision}, vol.~4, 2002, pp. 776--792.

\bibitem{ShalevShwartz2004}
S.~Shalev-Shwartz, Y.~Singer, and A.~Y. Ng, ``Online and batch learning of
  pseudo-metrics,'' in \emph{21st International Conference on Machine
  Learning}, 2004, p.~94.

\bibitem{Weinberger2005}
K.~Q. Weinberger, J.~Blitzer, and L.~K. Saul, ``Distance metric learning for
  large margin nearest neighbor classification,'' in \emph{Advances in neural
  information processing systems}, 2005, pp. 1473--1480.

\end{thebibliography}
	
\end{document}